# HEAVY TRAFFIC ANALYSIS OF OPEN PROCESSING NETWORKS WITH COMPLETE RESOURCE POOLING: ASYMPTOTIC OPTIMALITY OF DISCRETE REVIEW POLICIES

By Baris Ata and Sunil Kumar

*Northwestern University and Stanford University*

We consider a class of open stochastic processing networks, with feedback routing and overlapping server capabilities, in heavy traffic. The networks we consider satisfy the so-called complete resource pooling condition and therefore have one-dimensional approximating Brownian control problems. We propose a simple discrete review policy for controlling such networks. Assuming $2+\varepsilon$ moments on the interarrival times and processing times, we provide a conceptually simple proof of asymptotic optimality of the proposed policy.

**Contents**

**1. Introduction.** Stochastic processing networks have been extensively used to model manufacturing systems, computer systems, telecommunication networks and call-centers (see, e.g., [2, 9, 25, 40, 54]). One approach to designing control policies for such networks is the heavy traffic approximation approach pioneered by Harrison [14, 17]. This approach can be summarized by the following procedure; see [7, 15, 51]. (a) Formulate a stochastic network model and a notion of heavy traffic. (b) Formulate an approximating Brownian control problem for the network control problem, and reduce the dimension of this problem by deriving an equivalent workload formulation. (c) Analyze the Brownian control problem (or its equivalent workload formulation) and "interpret" its solution as a control policy for the original network. (d) Investigate the performance of the policy proposed in (c). In particular, determine whether it is asymptotically optimal in the heavy traffic limit.

Even though the heavy traffic approach has proved successful in many particular examples, [10, 21, 22, 29, 30, 31, 32, 37, 38, 46, 47, 48, 49], the complete procedure outlined above has not been resolved in general. The steps (a) and (b) have been resolved quite generally in the literature [7, 14, 17, 18, 20, 24, 32]. Steps (c) and (d) present several difficulties. First, the approximating Brownian control problem is not always analytically tractable. Second, even when the Brownian control problem is tractable, it is not easy to interpret its solution as a control policy for the original processing network. Finally, there are very few proofs of asymptotic optimality of the interpreted policy even when such an interpretation has been advanced [1, 16, 27, 28, 30, 35, 38, 42].

In this paper we will carry out all four steps (a)–(d) for a large class of open network models. The crucial assumption for our analysis is what Harrison and Lopez [19] called complete resource pooling (CRP), extended to a more general network setting similar to that considered by Bramson and Williams [8]. Roughly speaking, the CRP assumption requires enough overlap in the processing capabilities of the various servers to ensure that their capacities are exchangeable or transferable in the heavy traffic limit. Another way of saying this is that the equivalent workload formulation (EWF) of the approximating Brownian control problem referred to in (b) is one dimensional. The CRP assumption can be verified *a priori* by solving a linear program involving the first-order data of the stochastic network. We have linear holding cost function as the only economic element of our general model.

For networks that satisfy the CRP condition, the associated Brownian control problem is essentially trivial to solve. But, its solution, which prescribes keeping all but one buffer level at zero while keeping all servers busy, is not easy to interpret in the original stochastic network. We provide an interpretation based on the discrete review approach of Harrison [15] and Maglaras [33, 34]. Our interpretation provides a policy that reviews the contents of the buffers at discrete points in time, computes a processing plan based on the observed contents (interpreting zeros in the Brownian problem as small safety stocks), and implements this plan in open loop fashion. Moreover, even though the heavy traffic assumption is necessary for the formal analysis, the policy proposed can be adapted for the case where the system is near heavy traffic as well. This point will be elaborated on later.

All discrete review policies described to date in the literature of heavy traffic network control [15, 33, 34, 39] require the system manager to solve a new linear programming problem at each review point. However, by fully exploiting the special structure of CRP networks, we arrive at a far simpler type of policy, as follows. First, one solves a single linear program "off-line," this being essentially the static planning problem that one uses to verify the CRP assumption, and computes the inverse of a square matrix derived from its optimal basis matrix. Then, to determine the optimal processing plan at each review point, one simply multiplies a vector of observed buffer contents by that inverse matrix. Given the very simple structure of this policy, and fully exploiting the CRP assumption, we are able to prove asymptotic optimality in a conceptually simple fashion. Our proof uses the roadmap given by Bramson and Williams [5, 52, 53], but it does not invoke fluid limits. Rather, the proof establishes state space collapse (i.e., all buffer levels but one are zero in the diffusion limit) directly, and then uses the continuity of the one-dimensional reflection map to establish convergence to the desired limiting diffusion.

The class of network control problems whose equivalent workload formulation is one dimensional has received considerable research attention, see [1, 16, 19, 36, 44]. Two recent and major contributions to this area are the papers by Stolyar [44] and Mandelbaum and Stolyar [36]. These papers consider parallel server systems under the complete resource pooling assumption and establish the asymptotic optimality of the max-weight scheduling rule and the generalized $c\mu$ rule, respectively. Furthermore, the policies proposed in those papers do not require the knowledge of the external arrival rates. Our paper extends the body of knowledge in this area in several ways. First, we consider discretionary feedback routing which extends the network topology considered in [1, 19, 36]. Second, as mentioned above, the policy proposed is simple enough to be implementable on-line, since it only involves multiplication with a precomputed matrix at every review point. It is also simpler, if less elegant, than continuous review policies that require constant monitoring of the state; see [1, 35]. Third, we consider linear holding costs, which makes the translation in (c) harder since the optimal controls tend to achieve "corner solutions" in the state space rather than "interior solutions" as in [36, 44]. Therefore, one has to worry about interpreting zeros in the optimal Brownian prescriptions. Finally, our proof of asymptotic optimality is conceptually very simple, and only requires moments of order $2+\varepsilon$ for the interarrival and service times, rather than exponential moments as is usually assumed in the literature; see [1, 16, 35, 39].

**2. Description of the network model.** We assume that there are $p$ *servers* and $m$ *buffers*. The terms buffer and class will be used interchangeably. We assume that customers arrive to each buffer either from outside the network or from another buffer in the network. We assume that there are $n$ *activities*. Each activity is associated with a unique server and a unique buffer. When the server has expended sufficient time on an activity, a job either moves from the corresponding buffer to any of the $n$ buffers in a probabilistic manner or exits the system. Our model is a restricted version of the general processing network model of Harrison [17], as well as that of Bramson and Williams [8]. The stochastic assumptions will be made precise shortly. We let $s(j)$, $b(j)$ denote the server and the buffer associated with activity $j$, for $j=1,\dots,n$.

We describe the association between activities and resources by the capacity consumption matrix $A$, and the association between activities and buffers by the constituency matrix $C$. $A$ is a $p\times n$ matrix such that

$$A_{kj}=\begin{cases}1, & \text{if } s(j)=k,\\ 0, & \text{otherwise.}\end{cases} \tag{1}$$

$C$ is an $m\times n$ matrix such that

$$C_{ij}=\begin{cases}1, & \text{if } b(j)=i,\\ 0, & \text{otherwise.}\end{cases} \tag{2}$$

In our model $A,C$ are matrices of zeros and ones such that each contains exactly one nonzero entry in each column and at least one nonzero entry in each row. For an example of a network that fits in our modeling framework; see [26].

2.1. *Stochastic primitives.* Following the exposition of Bramson and Dai [6], we associate with each buffer $k=1,\dots,m$ a sequence of independent and identically distributed strictly positive random variables $\bar{u}_k=\{\bar{u}_k(i), i\ge 1\}$ and a $\lambda_k\ge 0$, where it is assumed that $\mathbf{E}(\bar{u}_k(1))=1$ for $k=1,\dots,m$. We allow $\lambda_k=0$ for some buffers but not all, and set $\mathcal{A}=\{k=1,\dots,m:\lambda_k\neq 0\}$.

We let $u_k(i) = \bar{u}_k(i)/\lambda_k$ for $k \in \mathcal{A}$ be the interarrival time between the $(i-1)$st and the $i$th externally arriving job at buffer $k$ for $k \in \mathcal{A}$, and $i = 1, 2, \ldots$ so that $\lambda_k$ is the external arrival rate to class $k$.

Similarly, for each activity $j = 1, \ldots, n$, we associate a sequence of strictly positive independent and identically distributed random variables $\bar{v}_j = \{\bar{v}_j(i), i \geq 1\}$ and a positive real number $m_k$. For each $i, j$ we let $v_j(i) = m_j \bar{v}_j(i)$ be the service time for the $i$th job processed by activity $j$. We also assume that $\mathbf{E}(\bar{v}_j(1)) = 1$ for $j = 1, \ldots, n$ so that for each activity $j$, $m_j$ is the mean service time of a job processed by activity $j$.

For each activity $j = 1, \ldots, n$, we also associate a sequence of independent and identically distributed $m$-dimensional random (routing) vectors $\phi_j = \{\phi_j(i), i \geq 1\}$. We let $\phi_j(i)$ be the routing vector of the $i$th job processed by activity $j$ and assume that it takes values in $\{e_0, e_1, \ldots, e_m\}$, where $e_0$ is the $m$-dimensional vector of zeros and $e_l$ is the $m$-dimensional vector with $l$th component 1 and other components 0. When $\phi_j(i) = e_0$, the $i$th job served by activity $j$ leaves the network, and when $\phi_j(i) = e_l$ for $l = 1, \ldots, m$, it next moves to buffer $l$. We let $P_{jl}$ denote the probability that $\phi_j(1) = e_l$ for $j = 1, \ldots, n$ and $l = 1, \ldots, m$ and we define the $n \times m$ activity-based routing matrix $P$ of our network as $P = (P_{jl})$.

We assume that the $2n + m$ random sequences, the stochastic increments,

$$\bar{u}_1, \ldots, \bar{u}_m, \qquad \bar{v}_1, \ldots, \bar{v}_n, \qquad \phi_1, \ldots, \phi_n, \tag{3}$$

are all defined on the same probability space and they are mutually independent. We specify the moment assumptions on the stochastic increments precisely as follows.

Moment assumptions. We assume that there exists an $\varepsilon_1 > 0$ such that

$$\mathbf{E}|\bar{u}_k(1)|^{2+2\varepsilon_1} < \infty \qquad \text{for } k = 1, \ldots, m, \tag{4}$$

and

$$\mathbf{E}|\bar{v}_j(1)|^{2+2\varepsilon_1} < \infty \qquad \text{for } j = 1, \ldots, n. \tag{5}$$

Our moment assumptions are weaker than the assumptions [16], Maglaras [33, 34] and Bell and Williams [1] used in analyzing asymptotically optimal policies.

We define the stochastic primitive processes of our network as the cumulative arrival, cumulative service and cumulative routing processes, which are, in turn, defined by the sums

$$U_k(l) = \sum_{i=1}^{l} u_k(i), \qquad V_j(l) = \sum_{i=1}^{l} v_j(i), \qquad \Phi_j(l) = \sum_{i=1}^{l} \phi_j(i), \tag{6}$$

where $j = 1, \ldots, n$, $k = 1, \ldots, m$ and $l = 1, 2, \ldots$.

For each activity $j$, $V_j(k)$ is the total amount of service required for the first $k$ jobs processed by that activity. We also define the renewal processes $E_k = \{E_k(t), t \geq 0\}$ and $S_j = \{S_j(t), t \geq 0\}$ associated with the external input process for class $k$ and service completion process for activity $j$, respectively, for $k \in \mathcal{A}$ and $j = 1, \ldots, n$ as follows:

$$E_k(t) = \max\left\{l \geq 0 : \sum_{i=1}^{l} u_k(i) \leq t\right\} \qquad \text{for } k \in \mathcal{A}, \tag{7}$$

$$S_j(t) = \max\left\{l \geq 0 : \sum_{i=1}^{l} v_j(i) \leq t\right\} \qquad \text{for } j = 1, \ldots, n. \tag{8}$$

We complete the definition of the first-order network data by introducing the input–output matrix $R$. First define $M$ as the diagonal matrix with entries $m_1, \ldots, m_n$, and $M^{-1}$ as its inverse. Next define the service rate for activity $j$ as $\mu_j = 1/m_j$ for $j = 1, \ldots, n$ and define the $m \times n$ input–output matrix as

$$R = (C - P')M^{-1}. \tag{9}$$

One can interpret the entry $R_{kj}$ as the average amount of material $k$ consumed by one unit of activity $j$, and negative values can be interpreted as production of material; see [14, 17].

2.2. *Assumptions on the first-order network data.* In the setting of conventional multiclass queueing networks, a queueing network is said to be in heavy traffic when all stations have utilization one. However, in the presence of dynamic routing decisions, the definition of heavy traffic is more subtle. Harrison [17] described the heavy traffic condition for such networks via a linear program, called static planning problem, involving the first-order network data. Along the same lines, we consider the static planning problem below:

STATIC PLANNING PROBLEM.

$$\text{minimize } \rho \tag{10}$$

$$\text{subject to } Rx = \lambda, Ax \leq \rho e, x \geq 0. \tag{11}$$

One can interpret $x_j$ as the long-run average rate at which activity $j$ is undertaken, expressed in units of activity per unit of time, and $\rho$ as a uniform upper bound on the utilization rates for the various resources under the processing plan $x$. In the static planning problem one seeks to minimize the maximal utilization rate $\rho$ subject to the requirement that average rates be nonnegative and that exogenously generated inputs be processed to completion without other inventories being generated.

Having introduced the static planning problem, we articulate the heavy traffic condition via static planning problem as follows:

HEAVY TRAFFIC ASSUMPTION. The static planning problem has a unique solution $(\rho^*, x^*)$. Moreover, that solution has $\rho^* = 1$ and $Ax^* = e$.

We will assume hereafter that the heavy traffic assumption is satisfied. One can interpret the heavy traffic assumption as follows: There is just one way of splitting arrivals in each input stream among available alternate routes so that no server is loaded beyond its capacity by the resultant flows, and each resource is critically loaded under this way of splitting the input streams.

In the solution of the static planning problem the activities for which $x_j^* > 0$ are called the basic activities and those for which $x_j^* = 0$ are called the nonbasic activities. We let $b$ denote the number of basic activities and relabel the activities so that activities $1, 2, \ldots, b$ are the basic ones. If there is degeneracy in the solution of the static planning problem, then the set of basic variables as defined above is not the same as the "basic" solution as understood in linear programming theory; see [3]. We partition $x^*$ as

$$x^* = \begin{bmatrix} x_B^* \\ x_N^* \end{bmatrix}, \tag{12}$$

where $x_B^*$ is the $b$-dimensional vector of nominal basic activity levels and $x_N^* = 0$. It will be also convenient for our later purposes to partition the input–output matrix $R$ and the capacity consumption matrix $A$ as follows:

$$R = [\,H \quad J\,] \quad \text{and} \quad A = [\,B \quad N\,], \tag{13}$$

where $H$ and $B$ both have $b$ columns and they correspond to the basic activities.

We also make the following natural assumption which simply says that each customer class is served by at least one basic activity.

ASSUMPTION BAB (Basic activity for each buffer). In the solution of the static planning problem, for each buffer $i$ there is at least one basic activity $j$ such that $H_{ij} > 0$.

In this paper, we only consider networks that satisfy the complete resource pooling (CRP) condition below.

COMPLETE RESOURCE POOLING ASSUMPTION. We assume that there is a full set of basic activities in the solution of the static planning problem. To be more specific, we assume that $p + m - b = 1$, where, as before, $p$ is the number of servers, $m$ is the number of buffers, and $b$ is the number of basic activities.

Analyzing a formal Brownian analog of the model in this section, Bramson and Williams [8] establish that (cf. Corollary (6.2) of [8]) the state descriptor is $d$-dimensional, where $d = p + m - b$. Thus, under the complete resource pooling assumption, the state descriptor is one dimensional and we exploit this fact significantly in our analysis.

We will also consider the dual linear program of the static planning problem which is given below.

Dual of the static planning problem.

$$\text{maximize } y'\lambda \tag{14}$$

$$\text{subject to } y'R \le \pi' A, \pi' e = 1 \quad \text{and} \quad \pi \ge 0. \tag{15}$$

Lemma 1. *There exist an $m$-dimensional vector $y$ and a $p$-dimensional vector $\pi$ such that $y$ and $\pi$ are both strictly positive and they satisfy the conditions below*:

$$y'H = \pi' B, \tag{16}$$

$$\pi' N \ge y' J. \tag{17}$$

This lemma is simply a restatement of Lemmas (7.2) and (7.6) of [8]. One can interpret $y_i$ as the workload contribution, or total work content, per class $i$ job and $\pi_k$ as the relative capacity of server $k$. Detailed interpretations of the vectors $y$ and $\pi$ are given in Section 4 of [17].

2.3. *Notation.* For each positive integer $k$, the $k$-dimensional Euclidean space will be denoted $R^k$; when $k = 1$, the superscript will be suppressed. $[x]$ denotes the integer part of a nonnegative real number $x$. Vectors will be column vectors unless indicated otherwise. Inequalities between vectors in $R^k$ should be interpreted componentwise. For $a, b \in R^k$, we shall use $a \vee b$ to denote the vector whose $i$th component is the maximum of $a_i$ and $b_i$ for $i = 1, \ldots, k$. Similarly, $a \wedge b$ will denote the componentwise minimum of $a$ and $b$. The superscript $'$ will be used to denote the operation of taking the transpose of a vector or matrix. For $x = (x_1, \ldots, x_k)' \in R^k$, we will use the norm $|x| = \max_{i=1}^k |x_i|$, and for the norm of an $k \times l$ matrix $A$, we will use $|A| = \max_{i=1}^k \sum_{j=1}^l |A_{ij}|$. For a vector $x \in R^k$, the $k \times k$ diagonal matrix whose diagonal entries are given by the components of $x$ will be denoted by $\operatorname{diag}(x)$. We define the "ball" around a point $z \in R^m$ of radius $a$ via

$$B_a(z) = \{q \in R^m : q_1 > z_1 - a, |q_i - z_i| < a \text{ for } i = 2, \ldots, m\}. \tag{18}$$

Let $(\Omega, \mathcal{F}, P)$ be a probability space. We denote filtrations on $(\Omega, \mathcal{F})$ by $\{\mathcal{F}_t, t \ge 0\}$. For each positive integer $k$, let $D^k$ be the space of all functions $\omega:[0,\infty)\to R^k$ that are right continuous on $[0,\infty)$ and have finite left limits on $(0,\infty)$. The identically zero function in $D^k$ will be denoted by $\mathbf{0}$. For $\omega\in D^k$ and $T\geq 0$, we let $\|\omega\|_T=\sup_{t\in[0,T]}|\omega(t)|$. Consider $D^k$ to be endowed with the usual Skorohod ($J1$) topology (see [4, 11]). Let $\mathcal{M}^k$ denote the Borel $\sigma$-algebra on $D^k$ associated with this topology. This is the same $\sigma$-algebra generated by the coordinate maps, that is, $\mathcal{M}^k=\sigma\{\omega(s):0\leq s<\infty\}$. Each continuous-time (stochastic) process in this paper will be assumed to be a measurable function from some probability space $(\Omega,\mathcal{F},P)$ into $(D^k,\mathcal{M}^k)$. A sequence of processes $\{W^i\}_{i=1}^{\infty}$ is said to be tight if and only if the probability measures induced by the sequence on $(D^k,\mathcal{M}^k)$ form a tight sequence. The notation "$W^i\Rightarrow W$" will mean that the probability measures induced by the $W^i$ on $(D^k,\mathcal{M}^k)$ converge weakly to the probability measure induced on $(D^k,\mathcal{M}^k)$ by $W$ as $i\to\infty$. For more on tightness and weak convergence of processes taking values in $D^k$ see [4, 11, 50].

**3. Scheduling controls and network dynamics.** We specify a scheduling policy or control by an $n$-dimensional continuous stochastic process $T=\{T(t),t\geq 0\}$, where $T_j(t)$ can be interpreted as the amount of service time devoted to activity $j$ by server $s(j)$ in $[0,t]$ for $j=1,\ldots,n$.

Having introduced the scheduling controls $T(\cdot)$, we define the performance related processes which will be driven by $T(\cdot)$. Let $Z_i(t)$ denote the number of class $i$ jobs in the network at time $t$ for $i=1,\ldots,m$; and $I_k(t)$ denote the cumulative idleness experienced by server $k$ up to time $t$ for $k=1,\ldots,p$. Also define the vector valued processes $Z,E,S$ associated with job-counts, external arrivals and service completions, respectively, in the obvious way. We also define cumulative idleness process as

$$I(t)=te-AT(t), \tag{19}$$

and the vector-valued deviation control process $Y(t)$ as

$$Y(t)=x^*t-T(t). \tag{20}$$

3.1. *Admissible policies and network dynamics.* In the literature of dynamic scheduling, most of the existing models deal only with work-conserving (or, nonidling) policies (see, e.g., [5, 27, 42, 52]). However, we will allow policies that may require the servers go idle even when there is work for them in the system. Indeed, the policy we will propose in the sequel will require idling the servers occasionally.

In this paper we restrict attention to so-called *head-of-line* policies. Loosely speaking, in these policies, when the server works on an activity, server effort is delivered solely to the job at the head of the line in the buffer corresponding to that activity. This assumption does preclude some policies such as processor sharing, but it is mathematically convenient. In particular, it allows us to describe the evolution of the buffer contents in the system by just specifying the cumulative time allocation process $T = \{T(t), t \geq 0\}$. For purposes of our analysis, it is sufficient to think of the head-of-line assumption as requiring $T$ to be admissible as defined below.

We call a scheduling policy $T = \{T(t), t \geq 0\}$ admissible if it maps stochastic primitive processes of our network $(E, S, \Phi_1, \ldots, \Phi_n)$ to $D^n[0, \infty)$, and the conditions below are satisfied. These conditions merely require that our model has reasonable system dynamics:

$$T(\cdot) \text{ is nondecreasing} \quad \text{and} \quad T(0) = 0, \tag{21}$$

$$A(T(t) - T(s)) \leq e(t - s) \qquad \text{for all } 0 \leq s \leq t < \infty, \tag{22}$$

$$Z(t) = Z(0) + E(t) + \sum_{j=1}^{n} \Phi_j(S_j(T_j(t))) - CS(T(t)), \tag{23}$$

$$Z(t) \geq 0 \qquad \text{for all } t \geq 0. \tag{24}$$

We note that the class of admissible policies is quite large. In particular, an admissible policy does not even have to be adapted to the natural filtration generated by the stochastic primitives. The conditions (21), (22) and (24) are quite natural in the sense that for any physical system, they have to be satisfied, and condition (22) can also be stated as $I(\cdot)$ is nondecreasing. The condition (23) reflects the fact that we only consider the head-of-the-line policies, and that we restrict attention to preemptive-resume policies.

3.2. *Objective.* Ideally, the objective of the system manager is to find the "optimal" scheduling controls. Assuming linear holding costs, a natural objective of such policy design is to find a nonanticipating control $T = \{T(t), t \geq 0\}$ that minimize the expected infinite horizon total discounted holding cost of the form

$$J_T = \mathbf{E}\left(\int_0^\infty e^{-\gamma s} h \cdot Z(s)\, ds\right), \tag{25}$$

where $Z(s)$ is the buffer content vector at time $s$, $h$ is the vector of holding cost rates, and $\gamma$ is the interest rate.

Another natural objective is to minimize the long-run average holding cost

$$J_T = \limsup_{t \to \infty} \mathbf{E}\left(\frac{1}{t} \int_0^t h \cdot Z(s)\, ds\right). \tag{26}$$

An even more ambitious objective is to find a nonanticipating policy that minimizes

$$P(h \cdot Z(s) > x) \qquad \text{for all } s > 0, x > 0 \tag{27}$$

among all admissible policies.

It is almost impossible in dynamic scheduling theory for stochastic processing networks to find solutions which can be described by a few parameters and are exactly optimal. However, lowering the aspirations in accordance with the general program laid out by Harrison [15] and Williams [7, 51], we relax the objective and seek "good" policies that are "asymptotically optimal." Namely, we construct a simple discrete review policy in the next section which turns out to be asymptotically optimal in a very strong sense. By asymptotically optimal, we mean that our policy achieves a lower bound on the system performance asymptotically; the precise definition of this statement and its proof will be provided in Section 5.

**4. Policy description.** In this section we describe our policy. To implement the policy, the system manager reviews the system status at discrete points in time and observes the buffer content levels. At each such review point, a nominal processing plan is derived for the ensuing period by a simple matrix computation and the resulting plan is implemented in open-loop fashion over the ensuing period. It is important to point out that we only use the first-order network data in specifying our policy. Even though performance of the policy depends critically on higher moments of the stochastic increments of the network, they are not used in describing the policy. Before describing our policy, we first present a policy-dependent performance bound in the next subsection to motivate our policy.

4.1. *A policy dependent performance bound.* We first define the workload process $W = \{W(t), t \geq 0\}$ as follows:

$$W(t) = \sum_{i=1}^{m} y_i Z_i(t) \qquad \text{for all } t \geq 0, \tag{28}$$

where $y$ is given by Lemma 1. Without loss of generality (by simply relabeling buffers), we can assume that

$$\frac{h_1}{y_1} \leq \frac{h_2}{y_2} \leq \cdots \leq \frac{h_m}{y_m}. \tag{29}$$

Clearly, we have the following lower bound on the instantaneous cost rate:

$$\sum_{i=1}^{m} h_i Z_i(t) \geq \frac{h_1}{y_1} W(t). \tag{30}$$

One can interpret the term on the right-hand side as the cost rate achieved by keeping all the workload in the "cheapest" buffer. Clearly, this lower bound is policy dependent because $W(t)$ depends on the policy employed. That is, we do not have a useful bound that works for all policies. However, we will provide an asymptotic bound that works for all policies (see Proposition 2).

4.2. *Description of the discrete review policy.* To motivate our choice of policy, we begin with an informal discussion of the policy in an idealized deterministic system in heavy traffic. The discrete review policy reviews the state of the system periodically. Consider the system at one such review point. Let $q$ denote the contents of the buffers at this review point. Based on $q$, we compute a processing plan $x$ to be implemented until the next review point, that is, over the next $l$ time units. Given the processing plan $x$, server $s(j)$ will spend $x_j l$ time units on activity $j$, serving buffer $b(j)$. We will choose $x$ so as to drive the state of the system at the end of the review period towards a target state $z$. Since the system is in heavy traffic, the workload cannot be decreased by any policy. So given (30), we would like to do two things with our choice of $x$ and $z$. First, we would like $Ax = e$ so as to ensure that the servers are fully utilized and thus prevent the workload from increasing. Second, we would like $z$ to be such that $z_2 = z_3 = \cdots = z_m = 0$. Of course, given an arbitrary $q$, it will not always be possible to achieve these objectives with a fixed period length $l$. But we would like to ensure that if the system started in a desirable state, that is, $q_2 = q_3 = \cdots = q_m = 0$, then it would continue to remain in a desirable state. In the idealized deterministic system, if $x$ is used as the processing plan for $l$ time units, the target state resulting is

$$z = q + \lambda l - Rxl,$$

where $R$ is given by (9). Let

$$\Pi = \begin{bmatrix} H & e_1 \\ B & 0 \end{bmatrix}, \tag{31}$$

where $H$ and $B$ are given by (13), and suppose for now that $\Pi^{-1}$ exists. If we used a processing plan $x = [x_B', x_N']'$, where $x_N = 0$ and a target state $z = (z_1, 0, 0, \ldots, 0)$ such that $x_B$ and $z_1$ are given by

$$\begin{bmatrix} x_B l \\ z_1 \end{bmatrix} = \Pi^{-1} \begin{bmatrix} q + \lambda l \\ el \end{bmatrix},$$

then we would have $z = q$ and $Ax = e$. That is, we would continue to stay in the desirable state and would have fully utilized the servers. In fact, the processing plan used would be the nominal processing plan $x^*$ from the solution of static planning problem [see (10) and (11)].

The prescription in the idealized system of keeping $q_2 = \cdots = q_m = 0$, while processing activities that drain these buffers according to the plan $x$, is not necessarily implementable in the original stochastic system. One way to adapt the policy is to specify a small safety stock $\theta_k$ for each buffer $k$ and to modify the processing plan and target state to

$$\begin{bmatrix} x_B l \\ z_1 - \theta_1 \end{bmatrix} = \Pi^{-1} \begin{bmatrix} q + \lambda l - \theta \\ el \end{bmatrix}. \tag{32}$$

This way, if we started in a desirable state $q_2 = \theta_2$, $q_3 = \theta_3, \ldots, q_m = \theta_m$ and $q_1 \geq \theta_1$, we would end up in a desirable state while fully utilizing the servers.

Finally, we need to specify about what to do if $q$ is far from a desirable state. Note that the processing plan in (32) tends to correct for deviations of $q$ from $\theta$. So if $q_2, \ldots, q_m$ were close to $\theta_2, \ldots, \theta_m$, we could still implement the plan from (32) and achieve target state $z$ that is closer to a desired state. However, the resulting processing plan $x$ would deviate from $x^*$, the solution of the static planning problem. Therefore, if $q$ were very far from $\theta$, the $x$ computed via (32) may be infeasible, that is, it may have negative components, capacity constraints may be violated, or may result in some components of $z$ being negative. In this case we need to stretch the length of the review period so that a feasible processing plan would still be able to achieve the desired target state. In fact, stretching the review period sufficiently long, we can achieve a desired state using only small perturbations of the nominal processing plan $x^*$.

We now make these heuristic ideas precise for the stochastic network model described in Section 2. Define constants

$$C_0 = \min_{1 \leq j \leq b, 1 \leq i \leq m} \frac{x_j^*}{|(\Pi^{-1}\binom{e_i}{0})_j|}, \tag{33}$$

$$C_1 = C_0 \left( \max_{j=1}^{m} \left\{ \frac{y_j}{y_1} \right\} \right) \tag{34}$$

and

$$\theta^* = n\left[ (2 + C_1 + C_0)\left( 1 \vee \max_{j=1}^{n} \{\mu_j\} \right) + 1 \right] e. \tag{35}$$

Given $l$, we choose the safety stock parameters $\theta$ as

$$\theta = \theta^* l, \tag{36}$$

and pick a perturbation constant $\delta > 0$ such that

$$\delta < \frac{C_0}{2m[1 + (\sum_{i=1}^m y_i)(\max_{k=1}^m \lambda_k)]}, \tag{37}$$

where $y$ is given by Lemma 1.

At the beginning of a review period, the system manager reviews the system status and observes $q$, the queue length vector, then determines the actual length of the planning period, $T^{exe}$ which may be different from the nominal length $l$, and the processing plan $x$ for the review period according to the prescription below. Then the servers start working in open-loop fashion by undertaking each (basic) activity for $x_j l$ time units or until the number of jobs processed using activity $j$ exceeds $q_{b(j)}/n$ in the review period. At the end of the review period, a new review period starts and the same procedure is repeated. We now describe the mechanics in a given review period.

We consider two cases: the first case corresponds to the case where the observed state and the target state are not too far apart, and the second to the case when they are.

CASE 1. $q \in B_{\delta l}(\theta)$. First, the system manager idles all the servers for $[y'(\theta - q)]^+$ to make sure that there is enough work in the system to be processed. All the processing activities are done in the next $l$ units of time. Therefore, we let the actual length of the review period be

$$T^{exe} = l + [y'(\theta - q)]^+. \tag{38}$$

We also set $x_N = 0$, $z_k = \theta_k$ for $k = 2, \ldots, m$ and determine the basic activity rates $x_B$ and target level for the first buffer $z_1$ by

$$\begin{bmatrix} x_B l \\ z_1 - \theta_1 \end{bmatrix} = \Pi^{-1} \begin{bmatrix} q + \lambda T^{exe} - \theta \\ el \end{bmatrix}. \tag{39}$$

Lemmas 2 and 3 establish that $\Pi$ is indeed invertible, and that $x$ and $z_1$ are well defined. Recall that we focus attention on preemptive-resume policies. The associated server $s(j)$ spends $x_j l$ units of time on activity $j$ during the review period, for $j = 1, \ldots, b$. Further we specify that for any basic activity $j$, the server $s(j)$ is not allowed to process more than $q_{b(j)}/n$ jobs from buffer $b(j)$ via activity $j$, where $n$ is the number of activities. This ensures that one activity does not overly drain a buffer and thus prevent other activities from being carried out. In case $q_{b(j)}/n$ jobs are processed by activity $j$ before the end of the review period, the server will simply go idle for the remaining time dedicated to activity $j$. By selecting $\theta$ large enough as in (35) and (36), we will see that this event will occur only with a small probability.

CASE 2. $q \notin B_{\delta l}(\theta)$. In this case the observed state $q$ is far enough from a desirable state so as to render the plan given by (39) infeasible, that is, $x$ does not satisfy $x \geq 0$, and $Ax \leq e$. We construct a feasible plan by first idling to accumulate work if necessary and then by "stretching" the length of the review period as follows. The system manager idles all the servers for $[y'(\theta - q)]^+$ to make sure that there is enough work in the system to be processed. We set $\tilde{q} = q + \lambda[y'(\theta - q)]^+$ and introduce stretching coefficient $C_s$, which will tell us by what factor we need to increase the length of the review period, given by

$$C_s = \max_{j=1,\ldots,b} \frac{|(1/l)[I, 0]\Pi^{-1}[\begin{smallmatrix} \theta - \tilde{q} \\ 0 \end{smallmatrix}])_j|}{x_j^*} \vee 1, \tag{40}$$

where $I$ is the $b \times b$ identity matrix. We let the actual length of the review period be

$$T^{exe} = [y'(\theta - q)]^+ + C_s l. \tag{41}$$

We also set $x_N = 0$ and $z_k = \theta_k$ for $k = 2, \ldots, m$. We then determine the target level for the first buffer and the basic activity rates, which will be undertaken for the last $C_s l$ time units of the review period, as

$$z_1 = \frac{y'(\tilde{q} - \theta)}{y_1} + \theta_1, \tag{42}$$

$$x_B = x_B^* \left(1 - \frac{1}{C_s}\right) + \frac{1}{C_s l}[I, 0]\Pi^{-1} \begin{bmatrix} \tilde{q} + \lambda l - \theta \\ el \end{bmatrix}, \tag{43}$$

where $I$ is the $b \times b$ identity matrix. Lemma 4 establishes that the processing plan just described is implementable in a deterministic setting.

As before, the policy now specifies that server $s(j)$ spends $C_s x_j l$ time units on activity $j$. If this results in buffer $b(j)$ emptying before the end of the review period, then the server simply idles for the remaining time devoted to activity $j$. This completes description of our discrete review policy. Case 2 is only needed for the sake of completeness, because it turns out that the probability that Case 2 ever arises under our policy vanishes in the heavy traffic limit. This observation, of course, simplifies our proofs significantly.

Having described the policy, we now state results showing that the policy is, indeed, well defined, and that it results in meaningful nominal allocations.

Lemma 2. *The policy matrix $\Pi$ is invertible.*

For the proof see Section A.1.

Lemma 3. *Given $q \in B_{\delta l}(\theta)$, $T^{exe} = l + [y'(\theta - q)]^+$ and if $x_B, z_1$ is given by*

$$\begin{bmatrix} x_B l \\ z_1 - \theta_1 \end{bmatrix} = \Pi^{-1} \begin{bmatrix} q + \lambda T^{exe} - \theta \\ el \end{bmatrix} \tag{44}$$

*and $z$ is given by*

$$z = q + \lambda T^{exe} - Rxl, \tag{45}$$

*then*

$$x_B \geq \tfrac{1}{2} x_B^* \tag{46}$$

*and*

$$z = \theta + \frac{1}{y_1} e_1 [y'(\theta - q)]^+. \tag{47}$$

For the proof see Section A.1.

LEMMA 4. *Given* $\tilde{q} \geq 0$ *such that* $y'\tilde{q} \geq y'\theta$, *and* $C_s, z_1, x_B$ *as in* (40), (42) *and* (43), *respectively, if* $z = \tilde{q} + \lambda C_s l - RxC_s l$, *where* $x = (x_B', 0)'$, *then*

$$x \geq 0, Ax \leq e \quad \text{and} \quad z_k = \theta_k \qquad \text{for } k \geq 2.$$

For the proof see Section A.1.

REMARK. It can be verified that the discrete review policy described above can still be implemented even if the heavy traffic assumption is not satisfied. To be more specific, suppose we perturb the arrival rates $\lambda$ to $\tilde{\lambda}$ and consider the static planning problem with $\tilde{\lambda}$. If this new static planning problem has a unique solution $(\tilde{\rho}, \tilde{x})$ with $\tilde{\rho} < 1$, and $A\tilde{x} = \tilde{\rho}e$, and if $\tilde{\lambda}$ is sufficiently close to $\lambda$ so that the same basis as in (10)–(13) is optimal, then the discrete review policy described above is still well defined for this case with the following minor change. One needs to replace the term "$el$" in (39) and (43) with the term "$\tilde{\rho}el$." Then, since Lemmas 2–4 only use the uniqueness assumption and not the assumption that $\rho = 1$, they can be modified to justify the validity of the policy. Furthermore, we can significantly improve the performance of this policy by using the excess capacity available. In particular, one could use the excess capacity to achieve complete state space collapse, where all the buffer levels are zero in the diffusion limit. We do not attempt a formal analysis of this case.

4.3. *Probability estimates for the discrete review policy.* In describing the discrete review policy, we only used the first-order network data. Even though description of the discrete review policy does not use anything but the first moments of the stochastic increments, its performance depends critically on the moments of higher order. As stated earlier in Section 2.1, we assume moments of order $2 + 2\varepsilon_1$ for some $\varepsilon_1 > 0$. In analyzing the performance of discrete review policies, Maglaras [33, 34] and Harrison [16] have imposed exponential moment assumptions on the stochastic increments, which is a strong assumption but it results in the tighter control of buffer content levels. On the other hand, we will need to take longer review periods because of the weaker moment assumptions, which will become apparent below. Our moment assumptions, along with the "long" review periods, simplifies the analysis.

In implementing the discrete review policy, system status is reviewed at discrete points in time, say $\tau_0, \tau_1, \tau_2, \ldots,$ where $\tau_0 = 0$ and the elements in this random sequence can be determined inductively.

We assume

$$Z(0) = \theta, \tag{48}$$

and that there are no partially completed jobs in the system at time zero and the arrival processes have no residual time at time 0. We define the set $\mathcal{N}_k$ for $k = 0, 1, 2, \ldots$ as

$$\mathcal{N}_k = A_k \cap B_k \cap C_k^c \cap D_k \cap E_k, \tag{49}$$

where

$$A_k = \{Z(\tau_k) \in B_{\delta l}(\theta), Z(\tau_{k+1}) \in B_{\delta l}(\theta)\}, \tag{50}$$

$$B_k = \left\{\max_{i=2}^{m} \sup_{\tau_k \le s \le \tau_{k+1}} Z_i(s) \le [2\theta^* + n(2|\mu| + 1) + 2|\lambda|(1 + y'\theta^*)]l\right\}, \tag{51}$$

$$\begin{aligned} C_k = {} & \{S_j(T_j(\tau_{k+1})) - S_j(T_j(\tau_k)) \ge (2\mu_j + 1)l \text{ for some } j\} \\ & \cap \{Z(\tau_k) \in B_{\delta l}(\theta)\}, \end{aligned} \tag{52}$$

$$D_k = \{\tilde{v}_j^{(k+1)} \le \delta\sqrt{l} \text{ for } j = 1, \ldots, n\} \cap \{\tilde{u}_i^{(k+1)} \le \delta\sqrt{l}\ \forall\, i \in \mathcal{A}\}, \tag{53}$$

$$E_k = \left\{\int_{\tau_k}^{\tau_{k+1}} \mathbb{1}_{\{W(s) > (y'\theta^* + mn|y|(2|\mu|+1) + 2|y|my'\theta^*|\lambda|)l\}}\, dI_W(s) = 0\right\}, \tag{54}$$

where $\tilde{v}_j^{(k)}$ is the residual service time for activity $j$ at time $\tau_k$ for $j = 1, \ldots, n$, and $k = 1, 2, \ldots$, and $\tilde{u}_i^{(k)}$ is the residual interarrival time for class $i$ at time $\tau_k$ for $i \in \mathcal{A}$, and $k = 1, 2, \ldots$, and

$$I_W(s) = \pi' I(s), \qquad s \ge 0. \tag{55}$$

$A_k$ tells us that we are in case 1 at the beginning and end of review period $k$, and $B_k$ specifies $Z_2$ through $Z_m$ have not grown too much during any part of the review period. If the number of jobs completed in a review period via activity $j$—where the initial queue length vector $q \in B_{\delta l}(\theta)$—exceeds $(2\mu_j + 1)l$, we will announce that to be a "coordination problem." That is, the event of coordination problems, $C_k$, in period $k$ is given by (52). This definition is, indeed, more stringent than necessary, because it might well be the case that for a sample path in this set, the servers are able to undertake their prescribed activity levels. However, it is quite cumbersome to enumerate all the possibilities which may lead to problems in undertaking the processing plan; and neither is it necessary. $D_k$ controls the overshoot of the residual interarrival and service times. Finally, $E_k$ helps us control the idleness incurred in a review period, and it is essential to observe for our future purposes that

$$\int_{\tau_k + [y'(\theta - Z(\tau_k))]^+}^{\tau_{k+1}} \mathbb{1}_{\{W(s) > (y'\theta^* + mn|y|(2|\mu|+1) + 2|y|my'\theta^*|\lambda|)l\}}\, dI_W(s) = 0,$$

on the set $C_k^c \cap \mathcal{N}_{k-1} \cap \cdots \cap \mathcal{N}_0$, which is a consequence of the fact that $I_W(\tau_k + [y'(\theta - Z(\tau_k))]^+) = I_W(\tau_{k+1})$ on that set. The latter assertion follows because the servers work continuously during $[\tau_k + [y'(\theta - Z(\tau_k))]^+, \tau_{k+1}]$ by our policy description. (As we restrict attention on the set $C_k^c \cap \mathcal{N}_{k-1} \cap \cdots \cap \mathcal{N}_0$, the servers will have enough input to work on during $[\tau_k + [y'(\theta - Z(\tau_k))]^+, \tau_{k+1}]$.)

PROPOSITION 1. *We fix an* $\varepsilon_1 \in (0,1)$ *such that* (4) *and* (5) *holds. Then*

$$\mathbf{P}(\mathcal{N}_k, \mathcal{N}_{k-1}, \dots, \mathcal{N}_0) \geq \left(1 - \frac{C}{l^{1+\varepsilon_1}}\right)^{k+1} \qquad \text{for all } l \geq C, k = 0, 1, 2, \dots, \tag{56}$$

*where* $C$ *is a constant independent of* $l$*; see* (138).

For the proof see Section A.3.

**5. Asymptotic analysis.** It is almost impossible in dynamic scheduling theory for stochastic processing networks to find solutions which can be described by a few parameters and are exactly optimal. However, lowering the aspirations in accordance with the general program laid out by Harrison [15] and Williams [7, 51], we relax the objective and seek "good" policies that are "asymptotically optimal" in the heavy traffic limit under diffusion scaling. In particular, we will establish the asymptotic optimality of the discrete review policy introduced earlier, provided its parameters are chosen correctly. The discrete review policy introduced in Section 4 will be denoted by $DR(l, \theta, \Pi)$, where $l$ is the (nominal) length of a review period, $\theta$ is the vector of safety stock levels, and $\Pi$ is the policy matrix. We will consider a sequence of systems indexed by the parameter $r$, and we will attach a superscript to note the dependence on $r$. The initial conditions and the parameters $l$ and $\theta$ will be varied with $r$ as below.

*Choice of policy parameters.* We fix an $\varepsilon_2 > 0$ such that $\varepsilon_2 < \varepsilon_1/3$, and choose the parameters $l, \theta$ of $DR(l, \theta, \Pi)$ for the $r$th system as follows:

$$l(r) = r^{1-\varepsilon_2}, \tag{57}$$

$$\theta(r) = \theta^* l(r). \tag{58}$$

*Initial conditions under scaling.* We assume for the $r$th system that

$$Z^r(0) = \theta(r), \tag{59}$$

and that there are no partially completed jobs in the system at time zero.

In order to analyze the asymptotic performance of the sequence of discrete review policies $\{DR(l(r), \theta(r), \Pi)\}_{r=1}^{\infty}$, we introduce the following diffusion scaled processes. Diffusion (or CLT) scaling is indicated by placing a hat over the process. We extend the definition of the scaled routing processes to all nonnegative times by making them piecewise constant. In defining the diffusion scaled quantities, we first center the processes, then accelerate the time by a factor of $r^2$ and normalize the space by a factor of $r$. A possible intuitive way to think about this type of scaling is to imagine that performance-relevant time spans are of order $r^2$ in the $r$th system and the natural units of measurement for queue lengths over such time spans are of order $r$.

*Diffusion scaled processes.*

$$\widehat{E}^r(t) = \frac{1}{r}[E(r^2t) - r^2\lambda t], \qquad t \ge 0, \tag{60}$$

$$\widehat{S}^r(t) = \frac{1}{r}[S(r^2t) - r^2\mu t], \qquad t \ge 0, \tag{61}$$

$$\widehat{\Phi}^r_j(t) = \frac{1}{r}[\Phi_j([r^2t]) - P'_j[r^2t]], \qquad t \ge 0, \tag{62}$$

$$\widehat{Z}^r(t) = \frac{1}{r}Z(r^2t), \qquad t \ge 0, \tag{63}$$

$$\widehat{W}^r(t) = \frac{1}{r}W(r^2t), \qquad t \ge 0, \tag{64}$$

$$\widehat{Y}^r(t) = \frac{1}{r}Y(r^2t), \qquad t \ge 0, \tag{65}$$

$$\widehat{I}^r(t) = \frac{1}{r}I(r^2t), \qquad t \ge 0, \tag{66}$$

where (62) defines the scaled routing vector for activity $j$, for $j = 1, \dots, n$.

Having introduced the diffusion scaled quantities, we now give a precise meaning to the term "asymptotic optimality."

*Asymptotic optimality.* A sequence of admissible policies $\{T^r_*(\cdot)\}_{r=1}^\infty$ is called asymptotically optimal if for any $t > 0$, $x > 0$,

$$\limsup_{r\to\infty} P(h \cdot \widehat{Z}^r_{T_*}(t) > x) \le \liminf_{r\to\infty} P(h \cdot \widehat{Z}^r_T(t) > x) \tag{67}$$

for any other sequence of admissible policies $\{T^r(\cdot)\}_{r=1}^\infty$.

As an aside for the reader familiar with heavy traffic literature, we note that the meaning of asymptotic optimality is different from, say, that considered by Bell and Williams [1]. In our setting the traffic intensity is always 1; it is with respect to policy parameters that we perform asymptotic analysis.

5.1. *Network dynamics under scaling.* To describe the evolution of system in terms of scaled processes, we need to introduce two additional scaled processes for which the time is accelerated by $r^2$ and the state space is

normalized by $r^2$:

$$\overline{T}^r(t) = \frac{1}{r^2} T(r^2 t), \qquad t \geq 0,$$

$$\overline{S}^r(t) = \frac{1}{r^2} S(r^2 t), \qquad t \geq 0.$$

It is now straightforward to derive

$$\widehat{Z}^r(t) = \widehat{X}^r(t) + R\widehat{Y}^r(t), \qquad t \geq 0, \tag{68}$$

$$\widehat{I}^r(t) = A\widehat{Y}^r(t), \qquad t \geq 0, \tag{69}$$

where

$$\widehat{X}^r(t) = \widehat{Z}^r(0) + \widehat{E}^r(t) + \sum_{j=1}^{n} \widehat{\Phi}_j^r(\overline{S}_j^r(\overline{T}_j^r(t))) - RM\widehat{S}^r(\overline{T}^r(t)), \qquad t \geq 0, \tag{70}$$

$$M = \operatorname{diag}(m). \tag{71}$$

We define $(n-b)$-dimensional vector $\eta$ such that

$$\eta' = \pi' N - y' J \tag{72}$$

and note that $\eta \geq 0$ by Lemma 1. Then observe that by premultiplying (68) by $y'$ and using (16), we have

$$\widehat{W}^r(t) = \widehat{X}_W^r(t) + \widehat{I}_W^r(t) + \widehat{I}_N^r(t), \qquad t \geq 0, \tag{73}$$

where

$$\widehat{X}_W^r(t) = y' \widehat{X}^r(t), \qquad t \geq 0, \tag{74}$$

$$\widehat{I}_W^r(t) = \pi' \widehat{I}^r(t), \qquad t \geq 0, \tag{75}$$

$$\widehat{I}_N^r(t) = -\eta' \widehat{Y}_N^r(t), \qquad t \geq 0. \tag{76}$$

5.2. *Convergence results.* In this section we present three convergence results regarding the scaled processes under $DR(l(r), \theta(r), \Pi)$ with $l(r), \theta(r)$ given by (57) and (58), respectively. These results are not only required to prove our main result [see Theorem 3], but also are interesting on their own right.

THEOREM 1 (State space collapse). *For any fixed time $T > 0$,*

$$\sup_{0 \leq s \leq T} \widehat{Z}_k^r(s) \Rightarrow \mathbf{0} \qquad \text{as } r \to \infty \text{ for } k = 2, \ldots, m$$

*under $\{DR(l(r), \theta(r), \Pi)\}_{r=1}^{\infty}$.*

Proof. We define the sequence of sets indexed by $r$ as

$$\mathcal{N}^r = \bigcap_{k=0}^{\lfloor r^2 T/l(r)\rfloor} \mathcal{N}_k^r.$$

We then have by Proposition 1 that

$$P(\mathcal{N}^r) \geq \left(1 - \frac{C}{l(r)^{1+\varepsilon_1}}\right)^{\lceil r^2 T/l(r)\rceil}.$$

Substituting $l(r) = r^{1-\varepsilon_2}$ [see (57)], gives

$$P(\mathcal{N}^r) \geq \left(1 - \frac{C}{(r^{(1-\varepsilon_2)})^{1+\varepsilon_1}}\right)^{\lceil (r^{1+\varepsilon_2})T\rceil}.$$

We note that

$$P(\mathcal{N}^r) \to 1 \qquad \text{as } r \to \infty,$$

because

$$\frac{C}{(r^{1-\varepsilon_2})^{1+\varepsilon_1}} T r^{1+\varepsilon_2} = \frac{TC}{r^{\varepsilon_1 - 2\varepsilon_2 - \varepsilon_1\varepsilon_2}} \to 0 \qquad \text{as } r \to \infty.$$

We then observe that for each fixed $r$, on the set $\mathcal{N}^r$, we have

$$\max_{k=2,\dots,m} \sup_{0\leq s\leq r^2 T} Z_k(s) \leq [2\theta^* + n(2|\mu| + 1) + 2|\lambda|(1 + y'\theta^*)]l(r).$$

Or, in terms of diffusion-scaled quantities, we have

$$\max_{k=2,\dots,m} \sup_{0\leq s\leq T} \widehat{Z}_k^r(s) \leq [2\theta^* + n(2|\mu| + 1) + 2|\lambda|(1 + y'\theta^*)]\frac{l(r)}{r}. \tag{77}$$

Clearly, for every $\varepsilon > 0$, we have

$$\begin{aligned} &P\left(\max_{k=2,\dots,m} \sup_{0\leq s\leq T} \widehat{Z}_k^r(s) > \varepsilon\right) \\ &\qquad \leq P((\mathcal{N}^r)^C) + P\left(\max_{k=2,\dots,m} \sup_{0\leq s\leq T} \widehat{Z}_k^r(s) > \varepsilon, \mathcal{N}^r\right). \end{aligned}$$

We can bound the right-hand side by using (77), which gives

$$\begin{aligned} &P\left(\max_{k=2,\dots,m} \sup_{0\leq s\leq T} \widehat{Z}_k^r(s) > \varepsilon\right) \\ &\qquad \leq P((\mathcal{N}^r)^C) + P\left(\frac{2\theta^* + n(2|\mu| + 1) + 2|\lambda|(1 + y'\theta^*)}{r^{\varepsilon_2}} > \varepsilon\right). \end{aligned}$$

Therefore, we have

$$P\left(\max_{k=2,\dots,m} \sup_{0\leq s\leq T} \widehat{Z}_k^r(s) > \varepsilon\right) \to 0 \qquad \text{as } r \to \infty.$$

Or, equivalently,

$$\sup_{0\le s\le T} \widehat{Z}_k^r(s) \Rightarrow 0 \qquad \text{as } r\to\infty \text{ for } k=2,\dots,m. \qquad \square$$

Before we state our next result, we first introduce the one-dimensional regulator map $(\psi,\varphi): D[0,\infty)\to D[0,\infty)$ by letting

$$\psi(x)(t) = -\inf_{0\le s\le t} x(s),$$

$$\varphi(x)(t) = x(t) + \psi(x)(t),$$

for all $x\in D[0,\infty)$.

THEOREM 2 (Convergence of scaled workload).

$$(\widehat{W}^r, \widehat{I}_W^r, \widehat{X}_W^r) \Rightarrow (W^*, I^*, X_W^*) \qquad \textit{as } r\to\infty$$

*under* $\{DR(l(r),\theta(r),\Pi)\}_{r=1}^{\infty}$, *where* $X_W^*$ *is a* $(0,\sigma^2)$ *Brownian motion starting at the origin,* $I^* = \psi(X_W^*)$, *and* $W^* = \varphi(X_W^*)$. *That is,* $W^*$ *is one-dimensional regulated Brownian motion. The variance parameter is* $\sigma^2 = y\Gamma y'$ *with* $\Gamma = \Gamma^0 + \sum_{j=1}^n x_j^*\Gamma^j$, *and* $\Gamma^0,\Gamma^1,\dots,\Gamma^n$ *are* $m\times m$ *covariance matrices defined as*

$$\Gamma_{jk}^0 = \lambda_k \operatorname{Var}(u_k(1))\mathbb{1}_{\{k=l\}},$$

$$\Gamma^j = \frac{1}{m_j}[\Omega^j + (R_jR_j')\operatorname{Var}(v_j(1))] \qquad \textit{for } j=1,\dots,n,$$

*where*

$$\Omega_{kl}^j = P_{jk}(\mathbb{1}_{\{k=l\}} - P_{jl}).$$

PROOF. We fix $T>0$ and analyze $DR(l(r),\theta(r),\Pi)$ over $[0,r^2T]$ for $r=1,2,\dots.$ As defined in the proof of Theorem 1, we let

$$\mathcal{N}^r = \bigcap_{k=0}^{\lfloor r^2T/l(r)\rfloor} \mathcal{N}_k^r,$$

and observe that for every sample path in the set $\mathcal{N}^r$, we have that

$$\int_0^{r^2T} \mathbb{1}_{\{W(s)>(y'\theta^*+mn|y|(2|\mu|+1)+2|y|my'\theta^*|\lambda|)l(r)\}}\, dI_W(s) = 0. \tag{78}$$

This follows immediately by construction of the sets $\mathcal{N}_k$, $k=0,1,\dots,\lfloor r^2T/l(r)\rfloor$. To be more specific, the set $E_k$ [see (54)], is constructed to make sure that we have (78) on the set $\mathcal{N}^r$.

One can equivalently represent (78) in terms of scaled quantities as follows:

$$\int_0^T \mathbb{1}_{\{\widehat{W}^r(s)>(y'\theta^*+mn|y|(2|\mu|+1)+2|y|my'\theta^*|\lambda|)l(r)/r\}}\,d\widehat{I}^r_W(s)=0. \tag{79}$$

Recall that [see (73)] under any admissible policy, we have that

$$\widehat{W}^r(t)=\widehat{X}^r_W(t)+\widehat{I}^r_W(t)+\widehat{I}^r_N(t).$$

Moreover, under $DR(l(r),\theta(r),\Pi)$, on the set $\mathcal{N}^r$, we have that

$$\widehat{I}^r_N=0 \qquad \text{for all } t\in[0,T].$$

This follows because $DR(l(r),\theta(r),\Pi)$ never uses the nonbasic activities. Consequently, under $DR(l(r),\theta(r),\Pi)$ we always have

$$\widehat{W}^r(t)=\widehat{X}^r_W(t)+\widehat{I}^r_W(t). \tag{80}$$

Having (79) and (80), we can invoke Lemma 7 which is stated and proved in Section A.2 to get, on the set $\mathcal{N}^r$,

$$\begin{aligned}\psi(\widehat{X}^r_W(s))&\le \widehat{I}^r_W(s)\\ &\le \psi(\widehat{X}^r_W(s))+\frac{(y'\theta^*+mn|y|(2|\mu|+1)+2|y|my'\theta^*|\lambda|)}{r^{\varepsilon_2}},\end{aligned} \tag{81}$$

$$\begin{aligned}\varphi(\widehat{X}^r_W(s))&\le \widehat{W}^r(s)\\ &\le \varphi(\widehat{X}^r_W(s))+\frac{(y'\theta^*+mn|y|(2|\mu|+1)+2|y|my'\theta^*|\lambda|)}{r^{\varepsilon_2}},\end{aligned} \tag{82}$$

for all $s\in[0,T]$.

The next step is to prove that $\overline{T}^r(\cdot)\Rightarrow x^*(\cdot)$ as $r\to\infty$ under $DR(l(r),\theta(r),\Pi)$, where $x^*(t)=x^*t$. To this end, we define the sequence processes $\overline{X}^r_W$, $\overline{W}^r$, $\overline{I}^r_W$, $\overline{X}^r$, $\overline{Z}^r$, $\overline{I}^r$, $\overline{Y}^r$ as follows:

$$\overline{X}^r_W(t)=\frac{1}{r}\widehat{X}^r(t), \qquad t\ge 0,$$

$$\overline{W}^r(t)=\frac{1}{r}\widehat{W}^r(t), \qquad t\ge 0,$$

$$\overline{I}^r_W(t)=\frac{1}{r}\widehat{I}^r_W(t), \qquad t\ge 0,$$

$$\overline{X}^r(t)=\frac{1}{r}\widehat{X}^r(t), \qquad t\ge 0,$$

$$\overline{Z}^r(t)=\frac{1}{r}\widehat{Z}^r(t), \qquad t\ge 0,$$

$$\overline{I}^r(t)=\frac{1}{r}\widehat{I}^r(t), \qquad t\ge 0,$$

$$\overline{Y}^r(t)=\frac{1}{r}\widehat{Y}^r(t), \qquad t\ge 0.$$

We note that $\overline{T}(s) \le s$, for all $t \in [0,T]$, by (22). We also recall [see (70)] that

$$\widehat{X}^r(t) = \widehat{Z}^r(0) + \widehat{E}^r(t) + \sum_{j=1}^n \widehat{\Phi}_j^r(\overline{S}_j^r(\overline{T}_j^r(t))) - RM\widehat{S}^r(\overline{T}^r(t)).$$

Using these, it is straightforward to derive

$$\|\widehat{X}^r\|_T \le \|\widehat{Z}^r(0)\| + \|\widehat{E}^r\|_T + \left\| \sum_{j=1}^n \widehat{\Phi}_j^r(\overline{S}_j^r) \right\|_{\|\overline{S}_j^r\|_T} + |RM| \|\widehat{S}^r\|_T.$$

It is also straightforward to show that the right-hand side converges in distribution to a nondegenerate limit. This can be proved by using the continuous mapping theorem, random time change theorem (cf. [4]) and the fact that $\overline{S}_j^r(s) \to \mu_j s$ as $r \to \infty$ for every $s \in [0,T]$ (cf. [23] for a proof). On the other hand, since $\overline{X}^r = \frac{1}{r}\widehat{X}^r$, we have that

$$\|\overline{X}^r\|_T \Rightarrow 0 \qquad \text{as } r \to \infty.$$

We also conclude by (74) that

$$\|\overline{X}_W^r\|_T \Rightarrow 0 \qquad \text{as } r \to \infty.$$

Since the one-dimensional regulator map commutes with scaling (see [50]), we immediately have, on the set $\mathcal{N}^r$, that

$$\psi(\overline{X}_W^r(s)) \le \overline{I}_W^r(s) \le \psi(\overline{X}_W^r(s)) + \frac{(y'\theta^* + mn|y|(2|\mu|+1) + 2|y|my'\theta^*|\lambda|)}{r^{\varepsilon_2}},$$

$$\varphi(\overline{X}_W^r(s)) \le \overline{W}^r(s) \le \varphi(\overline{X}_W^r(s)) + \frac{(y'\theta^* + mn|y|(2|\mu|+1) + 2|y|my'\theta^*|\lambda|)}{r^{\varepsilon_2}}$$

for all $s \in [0,T]$; and since $\psi, \varphi$ are continuous under Skorohod topology (cf. page 153 of [11]) and that $P(\mathcal{N}^r) \to 1$ as $r \to \infty$, we immediately have

$$\overline{I}_W^r \Rightarrow \mathbf{0} \qquad \text{as } r \to \infty,$$
$$\overline{W}^r \Rightarrow \mathbf{0} \qquad \text{as } r \to \infty.$$

Since we have $\pi > 0$, and $y > 0$ (cf. Lemma 1), we also have that

$$\overline{I}^r \Rightarrow \mathbf{0} \qquad \text{as } r \to \infty,$$
$$\overline{Z}^r \Rightarrow \mathbf{0} \qquad \text{as } r \to \infty.$$

We now prove that $\overline{Y}^r(\cdot) \Rightarrow \mathbf{0}$ as $r \to \infty$. Or, equivalently, $\overline{T}^r(\cdot) \Rightarrow x^*(\cdot)$ as $r \to \infty$. We first note that

$$Y = \begin{bmatrix} Y_B \\ Y_N \end{bmatrix} \quad \text{and} \quad T = \begin{bmatrix} T_B \\ T_N \end{bmatrix},$$

where $Y_N = T_N = \mathbf{0}$ under $DR(l(r),\theta(r),\Pi)$, because our policy never uses the nonbasic activities. Also, we have for any admissible policy that, by dividing both sides of (68) and (69) by $r$,

$$\overline{Z}^r(s) = \overline{X}^r(s) + R\overline{Y}^r(s),$$
$$\overline{I}^r(s) = A\overline{Y}^r(s).$$

Moreover, under $DR(l(r),\theta(r),\Pi)$ these reduce to the following:

$$\overline{Z}^r(s) = \overline{X}^r(s) + H\overline{Y}^r_B(s),$$
$$\overline{I}^r(s) = B\overline{Y}^r_B(s).$$

We can rewrite this in matrix notation as follows:

$$\begin{bmatrix} H \\ B \end{bmatrix} \overline{Y}^r_B(s) = \begin{bmatrix} \overline{Z}^r(s) - \overline{X}^r(s) \\ \overline{I}^r(s) \end{bmatrix}.$$

Or, we can write

$$\begin{bmatrix} H & -e_1 \\ B & 0 \end{bmatrix} \begin{bmatrix} \overline{Y}^r_B(s) \\ 0 \end{bmatrix} = \begin{bmatrix} \overline{Z}^r(s) - \overline{X}^r(s) \\ \overline{I}^r(s) \end{bmatrix}.$$

It is immediate from this and Lemma 2 that

$$\begin{bmatrix} \overline{Y}^r_B(s) \\ 0 \end{bmatrix} = \Pi^{-1} \begin{bmatrix} \overline{Z}^r(s) - \overline{X}^r(s) \\ \overline{I}^r(s) \end{bmatrix}.$$

Since $(\overline{Z}^r, \overline{X}^r, \overline{I}^r) \Rightarrow (\mathbf{0},\mathbf{0},\mathbf{0})$ as $r \to \infty$, we conclude by the continuous mapping theorem that $\overline{Y}^r_B(\cdot) \Rightarrow 0$ as $r \to \infty$, which in turn implies that

$$\overline{T}^r(\cdot) \Rightarrow x^*(\cdot) \qquad \text{as } r \to \infty.$$

We now immediately conclude by the representation (70) and the definition of $\widehat{X}^r_W$ [see (74)] and by the continuous mapping theorem, random time change theorem and functional central limit theorem for renewal processes (see [4]) that

$$\widehat{X}^r_W(\cdot) \Rightarrow X^*_W \qquad \text{as } r \to \infty,$$

where $X^*_W$ is a one-dimensional $(0,\, \sigma^2)$ Brownian motion. Deriving the expression for $\sigma^2$ is straightforward but tedious. It is outlined in [14] and so we will not repeat here.

We also conclude by the continuous mapping theorem that

$$\psi(\widehat{X}^r_W) \Rightarrow \psi(X^*_W) \qquad \text{as } r \to \infty,$$
$$\varphi(\widehat{X}^r_W) \Rightarrow \varphi(X^*_W) \qquad \text{as } r \to \infty,$$

because $\psi, \varphi$ are continuous under Skorohod topology; see [11].

Finally, since $P(\mathcal{N}^r) \to 1$ as $r \to \infty$ and (81), (82) holds on the set $\mathcal{N}^r$, we conclude by Skorohod representation theorem (see [4]) that

$$\widehat{I}_W^r \Rightarrow \psi(X_W^*) \qquad \text{as } r \to \infty,$$
$$\widehat{W}^r \Rightarrow \varphi(X_W^*) \qquad \text{as } r \to \infty. \qquad \square$$

COROLLARY 1 (Convergence of scaled queue lengths).

$$\widehat{Z}^r \Rightarrow Z^* \qquad \textit{as } r \to \infty$$

*under* $\{DR(l(r), \theta(r), \Pi)\}_{r=1}^{\infty}$, *where*

$$Z^*(t) = \left(\frac{W^*(t)}{y_1}, 0, \ldots, 0\right).$$

PROOF. By (28) we have

$$\widehat{Z}_1^r(s) = \frac{1}{y_1}\left[\widehat{W}^r(s) - \sum_{i=2}^{m} y_i \widehat{Z}_i^r(s)\right].$$

We also have by Theorem 1 that

$$\sum_{i=2}^{m} y_i \widehat{Z}_i^r \Rightarrow \mathbf{0} \qquad \text{as } r \to \infty,$$

and by Theorem 2 that

$$\widehat{W}^r(\cdot) \Rightarrow W^* \qquad \text{as } r \to \infty.$$

Therefore, the result follows by the convergence together lemma (see [4]). $\square$

5.3. *An asymptotic performance bound.* In this section we develop an asymptotic lower bound on the cost rate for admissible sequence of policies. The following result provides an asymptotic bound on the cost rate:

PROPOSITION 2. *Given an arbitrary sequence of admissible poli- cies* $\{T^r(\cdot)\}_{r=1}^{\infty}$, *for each* $t > 0$, $x > 0$, *one has that*

$$\liminf_{r\to\infty} P(h \cdot \widehat{Z}_T^r(t) > x) \geq P\left(\frac{h_1}{y_1} W^*(t) > x\right) = 2N\left(\frac{-xy_1}{h_1 \sigma \sqrt{t}}\right),$$

*where* $N(\cdot)$ *is the cumulative distribution function for a standard normal random variable.*

For the proof see Section A.4.

Therefore, the process $\frac{h_1}{y_1} W^*(\cdot)$ gives an asymptotic lower bound on the achievable cost rate.

5.4. *Asymptotic optimality of the discrete review policy.* Consider the sequence of discrete review policies $\{DR(l(r), \theta(r), \Pi)\}$ with $l(r) = r^{1-\varepsilon_2}$, $\theta(r) = \theta^* l(r)$. Let $T_* = \{T_*^r(\cdot)\}_{r=1}^{\infty}$ denote the sequence of cumulative time allocations under $\{DR(l(r), \theta(r), \Pi)\}_{r=1}^{\infty}$. Our main result can then be stated as follows.

THEOREM 3 (Asymptotic optimality of the discrete review policy). *For each $t > 0$ and $x > 0$, we have*

$$\lim_{r\to\infty} P(h \cdot \widehat{Z}^r_{T_*}(t) > x) = P\left(\frac{h_1}{y_1} W^*(t) > x\right).$$

*Therefore, $\{DR(l(r), \theta(r), \Pi)\}_{r=1}^{\infty}$ is asymptotically optimal.*

PROOF. By Corollary 1 we have

$$h \cdot \widehat{Z}^r_{T_*} \Rightarrow \frac{h_1}{y_1} W^* \qquad \text{as } r \to \infty.$$

Given Proposition 2, the result follows (see [4]). $\square$

## APPENDIX: TECHNICAL PROOFS

A.1. *Proofs in Section 4.2.* We now present the proofs of technical results presented in the text. We start by stating and proving a lemma which will be useful in proving Lemma 2.

LEMMA 5. *The matrix $\left[\begin{smallmatrix} H & 0 \\ B & -e \end{smallmatrix}\right]$ is invertible.*

PROOF. Suppose it is not invertible, then there exists $x_B \in R^b, \alpha \in R$ such that

$$\begin{bmatrix} H & 0 \\ B & -e \end{bmatrix} \begin{bmatrix} x_B \\ \alpha \end{bmatrix} = 0, \tag{83}$$

$$\begin{bmatrix} x_B \\ \alpha \end{bmatrix} \neq 0. \tag{84}$$

Clearly, $Hx_B = 0$, $Bx_B = \alpha e$ and using (16), we immediately conclude $\alpha = 0$. Therefore, we have

$$\begin{bmatrix} H & 0 \\ B & -e \end{bmatrix} \begin{bmatrix} x_B \\ 0 \end{bmatrix} = 0. \tag{85}$$

Recall that we denote the solution of the static planning problem by $x^*$, where

$$x^* = \begin{bmatrix} x_B^* \\ 0 \end{bmatrix},$$

and $x_B^* > 0$. Therefore, without loss of generality we can assume that

$$x_B^* + x_B \geq 0. \tag{86}$$

Also by using the heavy traffic assumption and (85), (86), we write

$$\begin{bmatrix} H & 0 \\ B & -e \end{bmatrix} \begin{bmatrix} x_B^* + x_B \\ 0 \end{bmatrix} = \begin{bmatrix} \lambda \\ e \end{bmatrix},$$

$$\begin{bmatrix} x_B^* + x_B \\ 0 \end{bmatrix} \geq 0.$$

This implies $[\begin{smallmatrix} x_B^* + x_B \\ 0 \end{smallmatrix}]$, $\rho = 1$ is an optimal solution to the static planning problem. However, by the heavy traffic assumption, the solution of static planning problem is unique. That is, $x_B = 0$, which contradicts (84). ☐

PROOF OF LEMMA 2. We consider solving the following equation:

$$\begin{bmatrix} H & 0 \\ B & -e \end{bmatrix} \begin{bmatrix} x \\ \rho \end{bmatrix} = \begin{bmatrix} e_1 \\ 0 \end{bmatrix}, \tag{87}$$

which has a unique solution by Lemma 5. By Cramér's rule (see page 233 of [45]), we conclude

$$\rho = \frac{\det[\begin{smallmatrix} H & e_1 \\ B & 0 \end{smallmatrix}]}{\det[\begin{smallmatrix} H & 0 \\ B & -e \end{smallmatrix}]}.$$

We can write (87) as

$$Hx = e_1, \tag{88}$$

$$Bx = e\rho. \tag{89}$$

Premultiplying (88) by $y'$ gives $y'Hx = y_1$, and also by (89), (16) we have that

$$y'Hx = \pi' Bx = \pi' e\rho = \rho.$$

Therefore, we have $\rho = y_1$ and we conclude $\rho > 0$, because $y > 0$ by Lemma 1. Clearly, by (87), $\rho > 0$ implies

$$\det \begin{bmatrix} H & e_1 \\ B & 0 \end{bmatrix} \neq 0,$$

which in turn implies $\Pi$ is invertible. ☐

The following lemma is needed to prove Lemma 3:

LEMMA 6. *Given $q \in B_{C_0/(2m)l}(\theta)$, if $x_B, z_1$ is given by*

$$\begin{bmatrix} x_B l \\ z_1 - \theta_1 \end{bmatrix} = \Pi^{-1} \begin{bmatrix} q + \lambda l - \theta \\ el \end{bmatrix}, \tag{90}$$

*then* $\left[\begin{smallmatrix} x_B \\ z_1 \end{smallmatrix}\right] > 0$. *More precisely,* $\left[\begin{smallmatrix} x_B \\ z_1 \end{smallmatrix}\right] \geq \frac{1}{2}\left[\begin{smallmatrix} x_B^* \\ \theta_1 \end{smallmatrix}\right]$.

PROOF. We let $\alpha, \tilde{q}$ be such that $\alpha \geq 0$,

$$q = \tilde{q} + \alpha e_1$$

and

$$|\tilde{q}_i - \theta_i| < \frac{C_0}{2m} l \qquad \text{for } i = 1, \ldots, m.$$

Recall that $\theta > C_1 le$. We first prove the following two facts:

$$\Pi^{-1} \begin{bmatrix} e_1 \\ 0 \end{bmatrix} = \begin{bmatrix} 0 \\ 1 \end{bmatrix}, \tag{91}$$

$$\left( \Pi^{-1} \begin{bmatrix} e_i \\ 0 \end{bmatrix} \right)_{b+1} = \frac{y_i}{y_1}. \tag{92}$$

Equation (91) follows from Lemma 2 and the fact that

$$\Pi \begin{bmatrix} 0 \\ 1 \end{bmatrix} = \begin{bmatrix} e_1 \\ 0 \end{bmatrix}.$$

To prove (92), we let $x, \beta$ be such that

$$\begin{bmatrix} x \\ \beta \end{bmatrix} = \Pi^{-1} \begin{bmatrix} e_i \\ 0 \end{bmatrix}.$$

We want to prove that $\beta = \frac{y_i}{y_1}$. It follows from the equation immediately above and Lemma 2 that

$$Bx = 0, \tag{93}$$

$$Hx + \beta e_1 = e_i. \tag{94}$$

By (16) and (93), premultiplying (94) by $y'$ gives

$$\beta y_1 = y_i.$$

Thus, $\beta = \frac{y_i}{y_1}$.

We now note that

$$\begin{bmatrix} x_B l \\ z_1 - \theta_1 \end{bmatrix} = \Pi^{-1} \left( \begin{bmatrix} \lambda l \\ el \end{bmatrix} + \begin{bmatrix} \tilde{q} - \theta \\ 0 \end{bmatrix} + \alpha \begin{bmatrix} e_1 \\ 0 \end{bmatrix} \right).$$

We then have, by the heavy traffic assumption, Lemma 2 and (91) that

$$\begin{bmatrix} x_B l \\ z_1 - \theta_1 \end{bmatrix} = \begin{bmatrix} x_B^* l \\ 0 \end{bmatrix} + \Pi^{-1} \begin{bmatrix} \tilde{q} - \theta \\ 0 \end{bmatrix} + \alpha \begin{bmatrix} 0 \\ 1 \end{bmatrix}.$$

Since $\alpha \geq 0$, we also have

$$\begin{bmatrix} x_B l \\ z_1 - \theta_1 \end{bmatrix} \geq \begin{bmatrix} x_B^* l \\ 0 \end{bmatrix} + \sum_{i=1}^{m} (\tilde{q}_i - \theta_i) \left( \Pi^{-1} \begin{bmatrix} e_i \\ 0 \end{bmatrix} \right).$$

We now rewrite the inequality for each component; we first consider

$$z_1 - \theta_1 \geq \sum_{i=1}^{m} (\tilde{q}_i - \theta_i) \left( \Pi^{-1} \begin{bmatrix} e_i \\ 0 \end{bmatrix} \right)_{b+1}.$$

We immediately conclude by (92) and choice of $\tilde{q}, \alpha$ that

$$z_1 - \theta_1 \geq -\sum_{i=1}^{m} |\tilde{q}_i - \theta_i| \left( \frac{y_i}{y_1} \right) \geq -\frac{C_0}{2m} l \sum_{i=1}^{m} \frac{y_i}{y_1} \geq -\frac{C_0}{2} l \max_{i=1}^{m} \left\{ \frac{y_i}{y_1} \right\} = -\frac{C_1}{2} l.$$

Therefore,

$$z_1 \geq \theta_1 - \frac{C_1}{2} l \geq \frac{\theta_1}{2} > 0.$$

The second inequality follows because $\theta > C_1 l e$.

Also, we have for $j = 1, \ldots, b$,

$$x_j l \geq x_j^* l + \sum_{i=1}^{m} (\tilde{q}_i - \theta_i) \left( \Pi^{-1} \begin{bmatrix} e_i \\ 0 \end{bmatrix} \right)_j.$$

Clearly, we also have

$$x_j l \geq x_j^* l - \sum_{i=1}^{m} |\tilde{q}_i - \theta_i| \left| \left( \Pi^{-1} \begin{bmatrix} e_i \\ 0 \end{bmatrix} \right)_j \right| \geq x_j^* l - m \frac{C_0}{2m} l \frac{x_j^*}{C_0} = \frac{x_j^*}{2} l > 0,$$

where the second inequality follows by the choice of $\tilde{q}$ and definition of $C_0$ [see (33)]. □

Proof of Lemma 3. We let $\tilde{q} = q + \lambda [y'(\theta - q)]^+$. To prove $x_B \geq \frac{1}{2} x_B^*$, it suffices by Lemma 6 to check that $\tilde{q} \in B_{C_0/(2m)l}(\theta)$. We check this componentwise, first consider $i = 2, \ldots, m$,

$$|\tilde{q}_i - \theta_i| \leq |q_i - \theta_i| + \lambda_i [y'(\theta - q)]^+ \leq \delta l + \left( \max_{k=1}^{m} \lambda_k \right) \sum_{i=1}^{m} y_i |\theta_i - q_i|.$$

Or, we can write

$$|\tilde{q}_i - \theta_i| \leq \delta l \left( 1 + \left( \max_{k=1}^{m} \lambda_k \right) \sum_{i=1}^{m} y_i \right) \leq \frac{C_0}{2m} l.$$

The last inequality follows by the choice of $\delta$ [cf. (37)]. To complete the proof of the first part of the lemma, we also establish a similar inequality for the first component of $\tilde{q}$, which indeed holds trivially:

$$\tilde{q}_1 = q_1 + \lambda_1 [y'(\theta - q)]^+ \geq q_1 \geq \theta_1 - \delta l,$$

where the last inequality follows from the fact that $q \in B_{\delta l}(\theta)$.

To prove $z = \theta + \frac{1}{y_1}[y'(q-\theta)]^+ e_1$, we first observe by (44) and Lemma 2 that

$$Hx_B l + e_1(z_1 - \theta_1) = q + \lambda T^{exe} - \theta, \tag{95}$$

$$Bx_B l = el. \tag{96}$$

We observe that premultiplying (95) by $y'$ gives

$$y'(Hx_B l + e_1(z_1 - \theta_1)) = y'[q + \lambda l + \lambda [y'(\theta - q)]^+ - \theta].$$

Using (16), we write

$$\pi' B x_B l + y_1(z_1 - \theta_1) = y'q + y'\lambda l + (y'\lambda)[y'(\theta - q)]^+ - y'\theta.$$

We further simplify this by using (96), the heavy traffic assumption and Lemma 1 to get

$$y_1(z_1 - \theta_1) = y'(q - \theta) + [y'(\theta - q)]^+.$$

Therefore, we have

$$z_1 = \theta_1 + \frac{[y'(q-\theta)]^+}{y_1}.$$

Finally, by (45) and (95) we have

$$z = q + \lambda T^{exe} - Rxl = q + \lambda T^{exe} - Hx_B l = \theta + (z_1 - \theta_1)e_1,$$

and we conclude the proof by observing that

$$z_k = \theta_k \qquad \text{for } k = 2, \dots, m. \qquad \square$$

Proof of Lemma 4. We first argue that $x \geq 0$. For $j = 1, \dots, b$, one can express $x_j$ as follows [from (43)]:

$$\begin{aligned}
x_j &= x_j^*\left(1 - \frac{1}{C_s}\right) + \frac{1}{C_s l}\left([I,0]\Pi^{-1}\begin{pmatrix}\tilde{q} + \lambda l - \theta \\ el\end{pmatrix}\right)_j \\
&= x_j^*\left(1 - \frac{1}{C_s}\right) + \frac{1}{C_s l}\left([I,0]\Pi^{-1}\begin{pmatrix}\lambda l \\ el\end{pmatrix} + [I,0]\Pi^{-1}\begin{pmatrix}\tilde{q} - \theta \\ 0\end{pmatrix}\right)_j \\
&= x_j^*\left(1 - \frac{1}{C_s}\right) + \frac{1}{C_s l}\left(x_B^* l + [I,0]\Pi^{-1}\begin{pmatrix}\tilde{q} - \theta \\ 0\end{pmatrix}\right)_j
\end{aligned}$$

$$= x_j^* - \frac{x_j^*}{C_s} + \frac{x_j^*}{C_s} + \frac{1}{C_s}\left(\frac{1}{l}[I,0]\Pi^{-1}\begin{pmatrix}\tilde{q}-\theta\\0\end{pmatrix}\right)_j$$

$$\geq x_j^* + \frac{1}{C_s}(-x_j^* C_s),$$

where the last inequality follows from the definition of $C_s$.

To complete the proof, we first observe that

$$C_s l x_B = [I,0]\Pi^{-1}\begin{pmatrix}\tilde{q}+\lambda C_s l-\theta\\C_s el\end{pmatrix}. \tag{97}$$

Also, one can show by using (92), (13) and the heavy traffic assumption that

$$z_1 - \theta_1 = \frac{[y'(\tilde{q}-\theta)]^+}{y_1} = [e'_{b+1}]\Pi^{-1}\begin{pmatrix}\tilde{q}+\lambda C_s l-\theta\\C_s el\end{pmatrix}, \tag{98}$$

where $e_{b+1}$ is a $(b+1)$-vector whose first $b$ entries are zeros and last entry is one. By combining (97) and (98), we arrive at the following:

$$\begin{bmatrix}x_B C_s l\\z_1-\theta_1\end{bmatrix} = \Pi^{-1}\begin{bmatrix}\tilde{q}+\lambda C_s l-\theta\\eC_s l\end{bmatrix},$$

which can equivalently be written as

$$\begin{bmatrix}H & e_1\\B & 0\end{bmatrix}\begin{bmatrix}x_B C_s l\\z_1-\theta_1\end{bmatrix} = \begin{bmatrix}\tilde{q}+\lambda C_s l-\theta\\eC_s l\end{bmatrix}.$$

One can further write this as follows:

$$H x_B C_s l + e_1(z_1-\theta_1) = \tilde{q}+\lambda C_s l-\theta, \tag{99}$$

$$B x_B C_s l = C_s l e. \tag{100}$$

Because $x_N = 0$, it is immediate from (100) that $Ax = e$. Finally, by using (99) and that $x_N = 0$, one has the following:

$$\begin{aligned}z &= \tilde{q}+\lambda C_s l - R x C_s l\\&= \tilde{q}+\lambda C_s l - H x_B C_s l\\&= \tilde{q}+\lambda C_s l - (\tilde{q}+\lambda C_s l-\theta-e_1(z_1-\theta_1))\\&= \theta + (z_1-\theta_1)e_1,\end{aligned}$$

which can also be written as $z_k = \theta_k$ for $k = 2,\ldots,m$. □

A.2. *An auxiliary lemma for Theorem 2.* We now state a technical lemma which is crucial for proving Theorem 2. The first part of this lemma is proved on pages 14 and 15 of [53] (cf. equation (8) of [53]), and the proof of the second part is very similar to that of first part; but we state and prove it—which is essentially the same as the proof given by Williams [53]—here for the sake of completeness.

LEMMA 7. *We let* $w, x, y \in D([0,T], R)$ *and* $\delta > 0$ *be such that:*

(i) $w(t) = x(t) + y(t)\ \forall t \in [0,T]$,
(ii) $w(t) \geq 0\ \forall t \in [0,T]$,

(iii) (a) $y(0) = 0$,
(b) $y$ *is nondecreasing,*
(c) $\int_{[0,T]} \mathbb{1}_{(\delta,\infty)}(w(t))\, dy(t) = 0$.

For $z \in R$, we let $z^- = \max\{0, -z\} = -\min\{0, z\}$. Then we define for all $t \in [0,T]$, that

$$\hat{y}(t) = \sup\{(x(s))^- : 0 \leq s \leq t\},$$
$$\tilde{y}(t) = \hat{y}(t) + \delta.$$

Then we have

$$y(t) \leq \tilde{y}(t) \qquad \forall t \in [0,T], \tag{101}$$

$$\hat{y}(t) \leq y(t) \qquad \forall t \in [0,T]. \tag{102}$$

REMARK ON CONDITION (iii)(c) OF LEMMA 7. One needs to be careful in interpreting the Lebesgue–Stieltjes integral, because the jumps in $x(\cdot)$ makes the interpretation quite subtle. We refer the interested reader to the remark on page 12 of [53] for a discussion of these issues.

PROOF OF LEMMA 7. We first present the proof of (101), which is, indeed, taken directly from [53]. We let $\varepsilon > 0$ and $\tau_\varepsilon = \inf\{t \in [0,T] : y(t) > \tilde{y}(t) + \varepsilon\}$ with $\inf \varnothing = \infty$. If $\tau_\varepsilon < \infty$, then $y(\tau_\varepsilon) \geq \tilde{y}(\tau_\varepsilon) + \varepsilon$ by the right continuity of paths, and

$$w(\tau_\varepsilon) = x(\tau_\varepsilon) + y(\tau_\varepsilon) \geq x(\tau_\varepsilon) + \tilde{y}(\tau_\varepsilon) + \varepsilon = x(\tau_\varepsilon) + \hat{y}(\tau_\varepsilon) + \delta + \varepsilon \geq \delta + \varepsilon.$$

However, by definition of $\tau_\varepsilon$ and (iii)(a) and (b), $y$ must either have a positive jump at time $\tau_\varepsilon$ or $y$ must have a point of increase to the right there. In either case, since $w(\tau_\varepsilon) > \delta$, this contradicts (iii)(c). Thus, $\tau_\varepsilon = \infty$ for each $\varepsilon > 0$, and, hence, $y(t) \leq \tilde{y}(t)\ \forall t \in [0,T]$ as desired.

We now prove (102). To this end, we first define $\hat{w}(t) = x(t) + \hat{y}(t)$ for $t \in [0,T]$. It is well known that (cf. [41]) $\hat{w}, x, \hat{y}$ jointly satisfy:

(i)′ $\hat{w}(t) = x(t) + \hat{y}(t)\ \forall t \in [0,T]$,
(ii)′ $\hat{w}(t) \geq 0\ \forall t \in [0,T]$,

(iii)′ (a) $\hat{y}(0) = 0$,
(b) $\hat{y}$ is nondecreasing,
(c) $\int_{[0,T]} \mathbb{1}_{(0,\infty)}(w(t))\, d\hat{y}(t) = 0$.

We let $\varepsilon > 0$, and $\tau_\varepsilon = \inf\{t \in [0,T] : \hat{y}(t) > y(t) + \varepsilon\}$ with $\inf \varnothing = \infty$. If $\tau_\varepsilon < \infty$, then $\hat{y}(\tau_\varepsilon) \geq y(\tau_\varepsilon) + \varepsilon$ by the right continuity of paths, and

$$\hat{w}(\tau_\varepsilon) = x(t) + \hat{y}(\tau_\varepsilon) \geq x(t) + y(\tau_\varepsilon) + \varepsilon = w(\tau_\varepsilon) + \varepsilon \geq \varepsilon > 0.$$

However, by definition of $\tau_\varepsilon$ and (iii)′(a) and (b), either $\hat{y}$ must have a positive jump at $\tau_\varepsilon$ or $\hat{y}$ must have a point of increase to the right there. In either case, since $\hat{w}(\tau_\varepsilon) > 0$, this contradicts (iii)′(c). Thus, $\tau_\varepsilon = \infty$ for each $\varepsilon > 0$ and, hence, $\hat{y}(t) \leq y(t)$ for all $t \in [0,T]$ as desired. □

A.3. *Proof of Proposition* 1. We now present the lemmas below which will be useful in proving Proposition 1.

LEMMA 8. *Given a sequence $\{Y_i\}_{i=1}^\infty$ of independent and identically distributed random variables with mean zero and*

$$E|Y_1|^{2+2\varepsilon_1} < \infty \qquad \text{for some } \varepsilon_1 > 0,$$

*we have*

$$E\left|\sum_{i=1}^N Y_i\right|^{2+2\varepsilon_1} \leq \left[\left(\frac{18(2+2\varepsilon_1)^{3/2}}{(1+2\varepsilon_1)^{1/2}}\right)^{2+2\varepsilon_1} E|Y_1|^{2+2\varepsilon_1}\right] N^{1+\varepsilon_1}.$$

PROOF. The result follows from Hölder's and Burkholder's inequalities. For a derivation of this result (indeed, a more general version) see equation (3.67) of [12]. □

We now prove a lemma regarding renewal processes associated with the external arrival and service processes.

LEMMA 9. *Given $\varepsilon > 0$ and $t > 2/\varepsilon$, we have*

$$P\Big(\sup_{0 \leq s \leq t} |E(s) - \lambda s| \geq \varepsilon t\Big) \leq \frac{C_2(\varepsilon)}{t^{1+\varepsilon_1}}, \tag{103}$$

$$P\Big(\sup_{x \geq 0, |x| \leq 1} |S(xt) - M^{-1}xt| \geq \varepsilon t\Big) \leq \frac{C_3(\varepsilon)}{t^{1+\varepsilon_1}}, \tag{104}$$

*where*

$$C_2(\varepsilon) = \frac{2+2\varepsilon_1}{2\varepsilon_1+1}\left[\frac{18(2+2\varepsilon_1)^{3/2}}{(1+2\varepsilon_1)^{1/2}}\right]^{2+2\varepsilon_1} \times \sum_{k \in \mathcal{A}} E|u_k(1) - Eu_k(1)|^{2+2\varepsilon_1}\left[\left(\frac{4\lambda_k^2(\lambda_k+\varepsilon)}{\varepsilon^2}\right)^{1+\varepsilon_1} + \left(\frac{4\lambda_k^3}{\varepsilon^2}\right)^{1+\varepsilon_1}\right],$$

$$C_3(\varepsilon) = \frac{2+2\varepsilon_1}{2\varepsilon_1+1}\left[\frac{18(2+2\varepsilon_1)^{3/2}}{(1+2\varepsilon_1)^{1/2}}\right]^{2+2\varepsilon_1}$$

$$\times \sum_{j=1}^{n} E|v_j(1) - Ev_j(1)|^{2+2\varepsilon_1} \left[ \left( \frac{4\mu_j^2(\mu_j + \varepsilon)}{\varepsilon^2} \right)^{1+\varepsilon_1} + \left( \frac{4\mu_j^3}{\varepsilon^2} \right)^{1+\varepsilon_1} \right].$$

PROOF. We will only present the proof of (103), because the proof of (104) is essentially the same. It is straightforward to arrive at the following:

$$P\Big(\sup_{0\le s\le t} |E(s) - \lambda s| \ge \varepsilon t\Big) \le \sum_{k\in\mathcal{A}} P\Big(\sup_{0\le s\le t} |E_k(s) - \lambda_k s| \ge \varepsilon t\Big).$$

Letting $\sigma_k = \inf\{s \in [0,t] : |E_k(s) - \lambda_k s| \ge \varepsilon t\}$, where $\inf \varnothing = \infty$, observe that

$$P\Big(\sup_{0\le s\le t} |E_k(s) - \lambda_k s| \ge \varepsilon t\Big) = P(\sigma_k \le t), \qquad k \in \mathcal{A}.$$

On the other hand, it is easy to see the following:

$$\begin{aligned}
\{|E_k(s) - \lambda_k s| \ge \varepsilon t\} &= \{E_k(s) - \lambda_k s \ge \varepsilon t\} \cup \{E_k(s) - \lambda_k s \le -\varepsilon t\} \\
&= \{E_k(s) \ge \lambda_k s + \varepsilon t\} \cup \{E_k(s) \le \lambda_k s - \varepsilon t\} \\
&\subset \{E_k(s) \ge \lfloor \lambda_k s + \varepsilon t \rfloor\} \cup \{E_k(s) \le \lceil \lambda_k s - \varepsilon t \rceil\} \\
&= \{U_k(\lfloor \lambda_k s + \varepsilon t \rfloor) \le s\} \cup \{U_k(\lceil \lambda_k s - \varepsilon t \rceil) \ge s\},
\end{aligned}$$

where $U(-k) = 0$ for $k \ge 0$.

Defining

$$\begin{aligned}
\sigma_k^{(1)} &= \inf\{s \in [0,t] : U_k(\lfloor \lambda_k s + \varepsilon t \rfloor) \le s\}, \qquad k \in \mathcal{A}, \\
\sigma_k^{(2)} &= \inf\{s \in [0,t] : U_k(\lceil \lambda_k s - \varepsilon t \rceil) \ge s\}, \qquad k \in \mathcal{A},
\end{aligned}$$

we can write $\sigma_k \ge \sigma_k^{(1)} \wedge \sigma_k^{(2)}$. Therefore,

$$P(\sigma_k \le t) \le P(\sigma_k^{(1)} \wedge \sigma_k^{(2)} \le t) \le P(\sigma_k^{(1)} \le t) + P(\sigma_k^{(2)} \le t).$$

We now analyze each of these terms separately. First, $\sigma_k^{(1)}$ can be re-expressed as follows:

$$\sigma_k^{(1)} = \inf\Big\{s \in [0,t] : U_k(\lfloor \lambda_k s + \varepsilon t \rfloor) - \frac{\lfloor \lambda_k s + \varepsilon t \rfloor}{\lambda_k} \le s - \frac{\lfloor \lambda_k s + \varepsilon t \rfloor}{\lambda_k}\Big\}.$$

Next, defining

$$\tilde{\sigma}_k^{(1)} = \inf\Big\{s \in [0,t] : U_k(\lfloor \lambda_k s + \varepsilon t \rfloor) - \frac{\lfloor \lambda_k s + \varepsilon t \rfloor}{\lambda_k} \le -\frac{\varepsilon t}{2\lambda_k}\Big\},$$

one has that $\tilde{\sigma}_k^{(1)} \le \sigma_k^{(1)}$. This follows because $\varepsilon t > 2$, and that in turn implies that $-\frac{\varepsilon t}{2\lambda_k} > s - \frac{\lfloor \lambda_k s + \varepsilon t \rfloor}{\lambda_k}$. Therefore,

$$P(\sigma_k^{(1)} \le t) \le P(\tilde{\sigma}_k^{(1)} \le t).$$

It follows that

$$P(\tilde{\sigma}_k^{(1)} \le t) \le P\biggl(\sup_{i=0,1,\dots,\lfloor \lambda_{k t+\varepsilon t}\rfloor} \biggl| U_k(i) - \frac{i}{\lambda_k}\biggr| \ge \frac{\varepsilon t}{2\lambda_k}\biggr).$$

It is straightforward to conclude by Markov's inequality (cf. page 39 of Ross [43]) that

$$P(\tilde{\sigma}_k^{(1)} \le t) \le \frac{E(\sup_{i=0,1,\dots,\lfloor \lambda_{k t+\varepsilon t}\rfloor} |U_k(i) - i/\lambda_k|)^{2+2\varepsilon_1}}{(\varepsilon t/(2\lambda_k))^{2+2\varepsilon_1}}.$$

We then use Doob's inequality (cf. page 15 of [12]) to get the following:

$$P(\tilde{\sigma}_k^{(1)} \le t) \le \frac{2+2\varepsilon_1}{1+2\varepsilon_1} \frac{E|U_k(\lfloor \lambda_k t + \varepsilon t\rfloor) - \lfloor \lambda_k t + \varepsilon t\rfloor/\lambda_k|^{2+2\varepsilon_1}}{(\varepsilon t/(2\lambda_k))^{2+2\varepsilon_1}}.$$

Finally, invoking Lemma 8 gives the following:

$$\begin{aligned} P(\tilde{\sigma}_k^{(1)} \le t) &\le \frac{2+2\varepsilon_1}{1+2\varepsilon_1}\biggl[\frac{18(2+2\varepsilon_1)^{3/2}}{(1+2\varepsilon_1)^{1/2}}\biggr]^{2+2\varepsilon_1} \\ &\quad \times E|u_k(1) - Eu_k(1)|^{2+2\varepsilon_1}\biggl[\frac{4(\lambda_k+\varepsilon)\lambda_k^2}{\varepsilon^2}\biggr]^{1+\varepsilon_1}\frac{1}{t^{1+\varepsilon_1}}. \end{aligned} \tag{105}$$

Similarly, we rewrite $\sigma_k^{(2)}$ as follows:

$$\sigma_k^{(2)} = \inf\biggl\{s \in [0,t] : U_k(\lceil \lambda_k s - \varepsilon t\rceil) - \frac{\lceil \lambda_k s - \varepsilon t\rceil}{\lambda_k} \ge s - \frac{\lceil \lambda_k s - \varepsilon t\rceil}{\lambda_k}\biggr\}.$$

Defining

$$\tilde{\sigma}_k^{(2)} = \inf\biggl\{s \in [0,t] : U_k(\lceil \lambda_k s - \varepsilon t\rceil) - \frac{\lceil \lambda_k s - \varepsilon t\rceil}{\lambda_k} \ge \frac{\varepsilon t}{2\lambda_k}\biggr\},$$

one has that $\tilde{\sigma}_k^{(2)} \le \sigma_k^{(2)}$. This follows because $\varepsilon t > 2$, and that, in turn, implies $\frac{\varepsilon t}{2\lambda_k} \le s - \frac{\lceil \lambda_k s - \varepsilon t\rceil}{\lambda_k}$. Therefore,

$$P(\sigma_k^{(2)} \le t) \le P(\tilde{\sigma}_k^{(2)} \le t).$$

It follows that

$$P(\tilde{\sigma}_k^{(2)} \le t) \le P\biggl(\sup_{i=0,1,\dots,\lceil \lambda_k t - \varepsilon t\rceil} \biggl| U_k(i) - \frac{i}{\lambda_k}\biggr| \ge \frac{\varepsilon t}{2\lambda_k}\biggr).$$

Then by Markov's inequality, one has that

$$P(\tilde{\sigma}_k^{(2)} \le t) \le \frac{E(\sup_{i=0,1,\dots,\lceil \lambda_k t - \varepsilon t\rceil} |U_k(i) - i/\lambda_k|)^{2+2\varepsilon_1}}{(\varepsilon t/(2\lambda_k))^{2+2\varepsilon_1}}.$$

One can further use Doob's inequality to arrive at the following:

$$P(\tilde{\sigma}_k^{(2)} \le t) \le \frac{2+2\varepsilon_1}{1+2\varepsilon_1} \frac{E|U_k(\lceil \lambda_k t - \varepsilon t\rceil) - \lceil \lambda_k t - \varepsilon t\rceil/\lambda_k|^{2+2\varepsilon_1}}{(\varepsilon t/(2\lambda_k))^{2+2\varepsilon_1}}.$$

Also, invoking Lemma 8 gives

$$\begin{aligned} P(\tilde{\sigma}_k^{(2)} \le t) &\le \frac{2+2\varepsilon_1}{1+2\varepsilon_1}\left[\frac{18(2+2\varepsilon_1)^{3/2}}{(1+2\varepsilon_1)^{1/2}}\right]^{2+2\varepsilon_1} \\ &\quad \times E|u_k(1) - Eu_k(1)|^{2+2\varepsilon_1}\left(\frac{4\lambda_k^3}{\varepsilon^2}\right)^{1+\varepsilon_1}\frac{1}{t^{1+\varepsilon_1}}. \end{aligned} \tag{106}$$

Finally, combining (105) and (106), for $k \in \mathcal{A}$, gives (103). □

The next lemma provides a similar probability estimate for our routing vectors, $\Phi_j(\cdot)$ for $j = 1, \ldots, n$.

LEMMA 10. *Given $\varepsilon > 0$ and for $j = 1, \ldots, n$, we have*

$$P\left(\sup_{i=0,1,\ldots,N} |\Phi_j(i) - P_j' i| \ge \varepsilon N\right) \le \frac{C_{4j}(\varepsilon)}{N^{1+\varepsilon_1}},$$

*where*

$$C_{4j}(\varepsilon) = \frac{2+2\varepsilon_1}{2\varepsilon_1+1}\left[\frac{18(2+2\varepsilon_1)^{3/2}}{(1+2\varepsilon_1)^{1/2}}\right]^{2+2\varepsilon_1} \sum_{k=1}^m E|\Phi_{jk}(1) - P_{jk}'|^{2+2\varepsilon_1}\left[\frac{1}{\varepsilon^{2+2\varepsilon_1}}\right].$$

PROOF. We first note that

$$P\left(\sup_{i=0,1,\ldots,N} |\Phi_j(i) - P_j' i| \ge \varepsilon N\right) \le \sum_{k=1}^m P\left(\sup_{i=0,1,\ldots,N} |\Phi_{jk}(i) - P_{jk}' i| \ge \varepsilon N\right).$$

By Markov's inequality (cf. page 39 of [43]), we write

$$P\left(\sup_{i=0,1,\ldots,N} |\Phi_j(i) - P_j' i| \ge \varepsilon N\right) \le \sum_{k=1}^m \frac{E(\sup_{i=0,1,\ldots,N}|\Phi_{jk}(i) - iP_{jk}'|)^{2+2\varepsilon_1}}{(\varepsilon N)^{2+2\varepsilon_1}}.$$

We then write by Doob's inequality (cf. page 15 of [12]),

$$P\left(\sup_{i=0,1,\ldots,N} |\Phi_j(i) - P_j' i| \ge \varepsilon N\right) \le \frac{2+2\varepsilon_1}{2\varepsilon_1+1} \sum_{k=1}^m \frac{E|\Phi_{jk}(N) - NP_{jk}'|^{2+2\varepsilon_1}}{N^{2+2\varepsilon_1}},$$

and the result follows by Lemma 8. □

PROOF OF PROPOSITION 1. We first review some essential notation. $\tau_0, \tau_1, \tau_2, \ldots$ are the review points with $\tau_0 = 0$. $\tau_k$ marks the beginning of the $k$th period for $k = 0, 1, 2, \ldots$; $x_j^{(i)}$ denotes rate of activity $j$ in period $i$ for $j=1,2,\ldots,n$, and $i=0,1,2,\ldots$. In particular, during review period $i$, server $s(j)$ (nominally) undertakes activity $j$ for $x_j^{(i)}l$ time units in Case 1, and $C_s x_j^{(i)} l$ time units in Case 2; and there are no partially completed jobs at time zero. Therefore, we have $x^{(0)}=x^*$, where $x^*$ is the solution of the static planning problem, and $\tau_1=l$.

We first prove that

$$P(\mathcal{N}_0)\geq 1-\frac{C}{l^{1+\varepsilon_1}}\qquad \text{for } l>C,$$

where $C$ is a constant, which depends on $\varepsilon$ but not on $l$, given by (138). Equivalently, we prove the following:

$$P(\mathcal{N}_0^c)\leq \frac{C}{l^{1+\varepsilon_1}}.$$

Since $\mathcal{N}_0^C=A_0^C\cup B_0^C\cup C_0\cup D_0^C\cup E_0^C$, we have

$$P(\mathcal{N}_0^C)\leq P(C_0)+P(D_0^C)+P(B_0^C\cap C_0^C)+P(A_0^C\cap C_0^C)+P(E_0^C\cap C_0^C).$$

We first consider $P(C_0)$:

$$\begin{aligned}P(C_0)&=P(S_j(x_j^*l)\geq(2\mu_j+1)l \text{ for some } j)\\&\leq P(S_j(l)\geq(2\mu_j+1)l \text{ for some } j)\\&\leq P(|S(el)-\mu l|\geq l).\end{aligned}$$

Therefore, by Lemma 9, we have that

$$P(C_0)\leq \frac{C_3(1)}{l^{1+\varepsilon_1}}\qquad \text{for } l>2. \tag{107}$$

We now consider $P(A_0^C\cap C_0^C)$. It follows that [recall that $Z(0)=\theta$]

$$\begin{aligned}Z(l)-\theta&=E(l)-\lambda l\\&\quad+\sum_{j=1}^n[\Phi_j(S_j(x_j^*l))-P_j'S_j(x_j^*l)]-RM[S(x^*l)-M^{-1}x^*l].\end{aligned}$$

Fix an $\varepsilon>0$ such that $\delta/6<\varepsilon<\delta/3$; and observe that

$$P(A_0^C\cap C_0^C)\leq p_I^0+p_{II}^0+p_{III}^0, \tag{108}$$

where

$$\begin{aligned}p_I^0&=P(|E(l)-\lambda l|\geq\varepsilon l),\\p_{II}^0&=P\left(\sum_{j=1}^n|\Phi_j(S_j(x_j^*l))-P_j'S_j(x_j^*l)|\geq\varepsilon l, S_j(x_j^*l)\leq(2\mu_j+1)l\ \forall j\right),\\p_{III}^0&=P(|RM||S(x^*l)-M^{-1}x^*l|\geq\varepsilon l).\end{aligned}$$

We will bound each of $p_I^0, p_{II}^0, p_{III}^0$ separately. We first observe by Lemma 9

$$p_I^0 = P(|E(l) - \lambda l| \geq \varepsilon l) \leq \frac{C_2(\varepsilon)}{l^{1+\varepsilon_1}} \qquad \text{for } l > \frac{12}{\delta}. \tag{109}$$

Then note that

$$\begin{aligned} p_{II}^0 &\leq \sum_{j=1}^n P\biggl(|\Phi_j(S_j(x_j^* l)) - P_j' S_j(x_j^* l)| \geq \frac{\varepsilon l}{n}, S_j(x_j^* l) \leq (2\mu_j + 1)l\biggr) \\ &\leq \sum_{j=1}^n P\biggl(\sup_{i=1,\dots,\lceil(2\mu_j+1)l\rceil} |\Phi_j(i) - iP_j'| \geq \frac{\varepsilon l}{n}\biggr) \\ &\leq \sum_{j=1}^n P\biggl(\sup_{i=1,\dots,\lceil(2\mu_j+1)l\rceil} |\Phi_j(i) - iP_j'| \geq \frac{\varepsilon}{n(2\mu_j+2)}\lceil(2\mu_j+1)l\rceil\biggr). \end{aligned}$$

Therefore, we conclude by Lemma 10 that

$$p_{II}^0 \leq \frac{1}{l^{1+\varepsilon_1}} \sum_{j=1}^n \frac{C_{4j}(\varepsilon/(n(2\mu_j+2)))}{(2\mu_j+1)^{1+\varepsilon_1}}. \tag{110}$$

Finally, we consider $p_{III}^0$:

$$p_{III}^0 \leq P\biggl(\sup_{\|\tilde{x}\| \leq 1, \tilde{x} \geq 0} |S(\tilde{x}l) - M^{-1}\tilde{x}l| \geq \frac{\varepsilon l}{|RM|}\biggr).$$

We then conclude by Lemma 9 that

$$p_{III}^0 \leq \frac{C_3(\varepsilon/|RM|)}{l^{1+\varepsilon_1}} \qquad \text{for } l > \frac{12|RM|}{\delta}. \tag{111}$$

Combining (108)–(111), we write

$$\begin{aligned} P(A_0^C \cap C_0^C) &\leq \biggl[C_2(\varepsilon) + C_3\biggl(\frac{\varepsilon}{|RM|}\biggr) + \sum_{j=1}^n \frac{C_4(\varepsilon/(n(2\mu_j+2)))}{(2\mu_j+1)^{1+\varepsilon_1}}\biggr]\frac{1}{l^{1+\varepsilon_1}}, \\ & \qquad l > \frac{12}{\delta}(1 \vee |RM|). \end{aligned} \tag{112}$$

We now consider $P(D_0^C)$ and note that

$$D_0^C = \bigcup_{j=1}^n \{\tilde{v}_j^{(1)} > \delta\sqrt{l}\,\} \cup \bigcup_{i \in \mathcal{A}} \{\tilde{u}_i^{(1)} > \delta\sqrt{l}\,\}.$$

Observe that for $i \in \mathcal{A}$, one has that

$$\begin{aligned} P(\tilde{u}_i^{(1)} > \delta\sqrt{l}\,) &= P(E_i(l + \delta\sqrt{l}\,) = E_i(l)) \\ &= \sum_{k=0}^\infty P(E_i(l + \delta\sqrt{l}\,) = k, E_i(l) = k) \end{aligned}$$

$$
\begin{aligned}
&= \sum_{k=0}^{\infty} P(U_i(k) < l, U_i(k+1) > l + \delta\sqrt{l}\,) \\
&\le \sum_{k=0}^{\infty} P(U_i(k) < l, u_i(k+1) > \delta\sqrt{l}\,) \\
&= \sum_{k=0}^{\infty} P(u_i(k+1) > \delta\sqrt{l}\,) P(E_i(l) = k) \\
&= P(u_i(1) > \delta\sqrt{l}\,) \sum_{k=0}^{\infty} P(E_i(l) = k) \\
&= P(u_i(1) > \delta\sqrt{l}\,).
\end{aligned}
$$

Therefore, we conclude by the Markov inequality that

$$
P(\tilde{u}_i^{(1)} > \delta\sqrt{l}\,) \le P(u_i(1) > \delta\sqrt{l}\,) \le \frac{E|u_i(1)|^{2+2\varepsilon_1}}{\delta^{2+2\varepsilon_1}} \frac{1}{l^{1+\varepsilon_1}}. \tag{113}
$$

Also, observe that for $j = 1, \ldots, n$, one has that

$$
\begin{aligned}
P(\tilde{v}_j^{(1)} > \delta\sqrt{l}\,) &= P(S_j(x_j^* l + \delta\sqrt{l}\,) = S_j(x_j^* l)) \\
&= \sum_{k=0}^{\infty} P(S_j(x_j^* l + \delta\sqrt{l}\,) = k, S_j(x_j^* l) = k) \\
&= \sum_{k=0}^{\infty} P(V_j(k) \le x_j^* l, V_j(k+1) > x_j^* l + \delta\sqrt{l}\,) \\
&\le \sum_{k=0}^{\infty} P(V_j(k) \le x_j^* l, v_j(k+1) > \delta\sqrt{l}\,) \\
&= \sum_{k=0}^{\infty} P(v_j(k+1) > \delta\sqrt{l}\,) P(S_j(x_j^* l) = k) \\
&= P(u_i(1) > \delta\sqrt{l}\,) \sum_{k=0}^{\infty} P(S_j(x_j^* l) = k) \\
&= P(v_j(1) > \delta\sqrt{l}\,).
\end{aligned}
$$

Therefore, we conclude by the Markov inequality that

$$
P(\tilde{v}_j^{(1)} > \delta\sqrt{l}\,) \le P(v_j(1) > \delta\sqrt{l}\,) \le \frac{E|v_j(1)|^{2+2\varepsilon_1}}{\delta^{2+2\varepsilon_1}} \frac{1}{l^{1+\varepsilon_1}}. \tag{114}
$$

We then conclude by (113) and (114) that

$$
P(D_0^C) \le \frac{1}{\delta^{2+2\varepsilon_1}} \left[ \sum_{j=1}^{n} E|v_j(1)|^{2+2\varepsilon_1} + \sum_{i \in \mathcal{A}} E|u_i(1)|^{2+2\varepsilon_1} \right] \frac{1}{l^{1+\varepsilon_1}}. \tag{115}
$$

Next, we consider $P(B_0^C \cap C_0^C)$. Since $Z(0) = \theta$, we have $\tau_1 = l$; and we have, for $0 \le s \le \tau_1$, that

$$Z(s) \le \theta + E(l) + \sum_{j=1}^{n} \Phi_j(S_j(x_j l)) \le \theta + E(l) + \left(\sum_{j=1}^{n} S_j(x_j l)\right) e.$$

Therefore, on the set $C_0^C$, we have, for $k = 2, \dots, m$, that

$$\sup_{0 \le s \le \tau_1} Z_k(s) \le \theta_k + E_k(l) + n(2\mu_j + 1)l.$$

Therefore, we note that

$$P(B_0^C \cap C_0^C) \le P(|E(l) - \lambda l| > |\lambda| l).$$

Hence, by Lemma 9, it is straightforward to conclude that

$$P(B_0^C \cap C_0^C) \le \frac{C_2(|\lambda|)}{l^{1+\varepsilon_1}} \qquad \text{for } l > \frac{2}{|\lambda|}. \tag{116}$$

Finally, we consider $P(E_0^C \cap C_0^C)$. However, by definition of the policy on the set $C_0^C$, the servers are never idled during $[\tau_0, \tau_1]$, and, hence, $I_W(l) = I_W(0) = 0$ on that set. The last assertion follows because the servers work continuously during $[0, l]$ by our policy description. (As we restrict attention on the set $C_0^C$, the servers will have enough input to work on during $[0, l]$.) Hence, trivially,

$$\int_{\tau_k}^{\tau_{k+1}} \mathbb{1}_{\{W(s) > (y'\theta^* + mn|y|(2|\mu|+1) + 2|y|my'\theta^*|\lambda|)l\}} \, dI_W(s) = 0 \qquad \text{on } C_0^C.$$

Therefore,

$$P(E_0^C \cap C_0^C) = 0. \tag{117}$$

Combining (107), (112), (115)–(117), we conclude that

$$P(\mathcal{N}_0^C) \le \frac{C}{l^{1+\varepsilon_1}} \qquad \text{for } l > 2 \vee \frac{12}{\delta} \vee \frac{12|RM|}{\delta} \vee \frac{2}{|\lambda|},$$

where the constant

$$C \ge C_3(1) + \left[C_2(\varepsilon) + C_3\left(\frac{\varepsilon}{|RM|}\right) + \sum_{j=1}^{n} \frac{C_4(\varepsilon/(n(2\mu_j + 2)))}{(2\mu_j + 1)^{1+\varepsilon_1}}\right]$$
$$+ \left[\sum_{j=1}^{n} E|v_j(1)|^{2+2\varepsilon_1} + \sum_{i \in \mathcal{A}} E|u_i(1)|^{2+2\varepsilon_1}\right] + C_2(|\lambda|)$$

[see (138)]. This provides the induction basis for proving Proposition 1. We now assume that (56) holds for $i = 0, 1, \dots, k-1$, and prove that it holds for $i = k$ to complete the proof. It suffices to prove

$$P(\mathcal{N}_k^c, \mathcal{N}_{k-1}, \dots, \mathcal{N}_0) \le \frac{C}{l^{1+\varepsilon_1}} P(\mathcal{N}_{k-1}, \dots, \mathcal{N}_0), \tag{118}$$

because

$$P(\mathcal{N}_k, \mathcal{N}_{k-1}, \ldots, \mathcal{N}_0) = P(\mathcal{N}_{k-1}, \ldots, \mathcal{N}_0) - P(\mathcal{N}_k^c, \mathcal{N}_{k-1}, \ldots, \mathcal{N}_0).$$

Therefore, we will show (118). Clearly, $\mathcal{N}_k^c = A_k^c \cup B_k^c \cup C_k \cup D_k^c \cup E_k^C$, and we have that

$$\begin{aligned} P(\mathcal{N}_k^C, \mathcal{N}_{k-1}, \ldots, \mathcal{N}_0) \le P(D_k^c, \mathcal{N}_{k-1}, \ldots, \mathcal{N}_0) &+ P(C_k, \mathcal{N}_{k-1}, \ldots, \mathcal{N}_0) \\ &+ P(B_k^c, C_k^c, \mathcal{N}_{k-1}, \ldots, \mathcal{N}_0) \\ &+ P(A_k^c, C_k^c, \mathcal{N}_{k-1}, \ldots, \mathcal{N}_0) \\ &+ P(E_k^c, C_k^c, \mathcal{N}_{k-1}, \ldots, \mathcal{N}_0). \end{aligned}$$

On the set $\mathcal{N}_0 \cap \mathcal{N}_1 \cap \cdots \cap \mathcal{N}_{k-1}$, we have $Z(\tau_i) \in B_{\delta l}(\theta)$ for $i = 1, \ldots, k$. Therefore, in each of these periods the system manager implements the Case 1 of $DR(l(r), \theta(r), \Pi)$. In particular,

$$T^{exe}(i) = \tau_i - \tau_{i-1} = l + [y'(\theta - Z(\tau_i))]^+ \qquad \text{for } i = 1, 2, \ldots, k, \tag{119}$$

$$\begin{bmatrix} x_B^{(i)} l \\ z_1^{(i)} - \theta_1 \end{bmatrix} = \Pi^{-1} \begin{bmatrix} Z(\tau_i) + \lambda T^{exe}(i) - \theta \\ el \end{bmatrix}, \tag{120}$$

and $x_N^{(i)} = 0$, and $x_B^{(i)} \ge (1/2) x_B^*$ by Lemma 3.

Therefore, to calculate the probability of an event $F \cap \mathcal{N}_0 \cap \mathcal{N}_1 \cap \cdots \cap \mathcal{N}_{k-1}$ for $i = 1, \ldots, k$, we first further partition this event with respect to the values of $Z(\tau_i) = q_i \in B_{\delta l}(\theta)$, $i = 1, \ldots, k$. That is, we fix $Z(\tau_i) = q_i \in B_{\delta l}(\theta)$ for $i = 1, \ldots, k$. It is important to observe that once we fix $Z(\tau_i) = q_i \in B_{\delta l}(\theta)$ for $i = 1, \ldots, k$, we also fix $\tau_1, \ldots, \tau_k$ and $x^{(1)}, \ldots, x^{(k)}$ on the set $\mathcal{N}_0 \cap \mathcal{N}_1 \cap \cdots \cap \mathcal{N}_{k-1}$ by (119) and (120). In other words, once we partition the event in question further with respect to values taken by $Z(\tau_i)$ for $i = 1, \ldots, k$, $\tau_i$ and $x^{(i)}$, for $i = 1, \ldots, k$, are determined by (119) and (120), and their values are deterministic (but different on each partition set) on each of these partition events whose union constitute the event in question. We then have

$$\begin{aligned} &P(F, \mathcal{N}_{k-1}, \ldots, \mathcal{N}_0) \\ &\quad = \sum_{(q_1, \ldots, q_k) \in B_{\delta l}(\theta) \times \cdots \times B_{\delta l}(\theta) \cap Z^{m \times k}} P(F, \mathcal{N}_{k-1}, \ldots, \mathcal{N}_0, \\ &\qquad\qquad\qquad\qquad\qquad\qquad Z(\tau_i) = q_i, i = 1, \ldots, k). \end{aligned}$$

Therefore, to simplify the notational burden, we will just pretend that $\tau_k$, $x^{(i)}$ are deterministic in the derivations below. However, we implicitly do the partitioning above to justify the calculations. We will illustrate this approach in estimating $P(D_k^c, \mathcal{N}_{k-1}, \ldots, \mathcal{N}_0)$, but this detailed approach will not be repeated in the interest of brevity. Nevertheless, one needs to go through the same steps to rigorously justify the calculations we will present. We first consider $P(D_k^c, \mathcal{N}_{k-1}, \dots, \mathcal{N}_0)$:

$$P(D_k^c, \mathcal{N}_{k-1}, \dots, \mathcal{N}_0) \le \sum_{j=1}^{n} P(\tilde{v}_j^{(k+1)} > \delta\sqrt{l}, \mathcal{N}_0, \dots, \mathcal{N}_{k-1}) + \sum_{i\in\mathcal{A}} P(\tilde{u}_i^{(k+1)} > \delta\sqrt{l}, \mathcal{N}_0, \dots, \mathcal{N}_{k-1}),$$

where

$$\begin{aligned}
&P(\tilde{u}_i^{(k+1)} > \delta\sqrt{l}, \mathcal{N}_0, \dots, \mathcal{N}_{k-1}) \\
&\quad = P(E_i(\tau_{k+1} + \delta\sqrt{l}\,) = E_i(\tau_{k+1}), \mathcal{N}_0, \dots, \mathcal{N}_{k-1}) \\
&\quad = \sum_{(q_1,\dots,q_k)\in B_{\delta l}(\theta)\times\cdots\times B_{\delta l}(\theta)} P(E_i(\tau_{k+1} + \delta\sqrt{l}\,) = E_i(\tau_{k+1}), \\
&\qquad\qquad \mathcal{N}_0, \dots, \mathcal{N}_{k-1}, Z(\tau_j) = q_j, j = 1, \dots, k) \\
&\quad = \sum_{(q_1,\dots,q_k)} P\Biggl(E_i\Biggl((k+1)l + \sum_{j=1}^{k}[y'(\theta - q_j)]^+ + \delta\sqrt{l}\Biggr) \\
&\qquad\qquad = E_i\Biggl((k+1)l + \sum_{j=1}^{k}[y'(\theta - q_j)]^+\Biggr) \\
&\qquad\qquad \mathcal{N}_0, \dots, \mathcal{N}_{k-1}, Z(\tau_j) = q_j, j = 1, \dots, k\Biggr) \\
&\quad \le \sum_{(q_1,\dots,q_k)} \sum_{\alpha=0}^{\infty} P(u_i(\alpha+1) > \delta\sqrt{l}\,) \\
&\qquad\qquad \times P\Biggl[E_i\Biggl((k+1)l + \sum_{j=1}^{k}[y'(\theta - q_j)]^+\Biggr) = \alpha, \\
&\qquad\qquad \mathcal{N}_0, \dots, \mathcal{N}_{k-1}, Z(\tau_j) = q_j, j = 1, \dots, k\Biggr] \\
&\quad = P(u_i(1) > \delta\sqrt{l}\,) \\
&\qquad \times \sum_{(q_1,\dots,q_k)} \sum_{\alpha=0}^{\infty} P(E_i(\tau_{k+1}) = \alpha, \mathcal{N}_0, \dots, \mathcal{N}_{k-1}, Z(\tau_j) = q_j, j = 1, \dots, k).
\end{aligned}$$

Therefore, we have that

$$P(\tilde{u}_i^{(k+1)} > \delta\sqrt{l}, \mathcal{N}_0, \ldots, \mathcal{N}_{k-1}) \leq P(u_i(1) > \delta\sqrt{l}\,) \times P(\mathcal{N}_0, \ldots, \mathcal{N}_{k-1}) \leq \frac{E|u_i(1)|^{2+2\varepsilon_1}}{\delta^{2+2\varepsilon_1}} \frac{1}{l^{1+\varepsilon_1}} P(\mathcal{N}_0, \ldots, \mathcal{N}_{k-1}).$$

Similarly, we write

$$\begin{aligned} &P(\tilde{v}_j^{(k+1)} > \delta\sqrt{l}, \mathcal{N}_0, \mathcal{N}_1, \ldots, \mathcal{N}_{k-1}) \\ &\quad = P\left(S_j\left(\sum_{i=0}^{k} x_j^{(i)} l + \delta\sqrt{l}\right) = S_j\left(\sum_{i=0}^{k} x_j^{(i)} l\right), \mathcal{N}_0, \mathcal{N}_1, \ldots, \mathcal{N}_{k-1}\right) \\ &\quad = \sum_{\alpha=0}^{\infty} P\left(S_j\left(\sum_{i=0}^{k} x_j^{(i)} l + \delta\sqrt{l}\right) = \alpha, S_j\left(\sum_{i=0}^{k} x_j^{(i)} l\right) = \alpha, \mathcal{N}_0, \mathcal{N}_1, \ldots, \mathcal{N}_{k-1}\right) \\ &\quad \leq \sum_{\alpha=0}^{\infty} P(v_j(\alpha+1) > \delta\sqrt{l}\,) \times P\left(S_j\left(\sum_{i=0}^{k} x_j^{(i)} l\right) = \alpha, \mathcal{N}_0, \mathcal{N}_1, \ldots, \mathcal{N}_{k-1}\right) \\ &\quad = P(v_j(1) > \delta\sqrt{l}\,) \sum_{\alpha=0}^{\infty} P\left(S_j\left(\sum_{i=0}^{k} x_j^{(i)} l\right) = \alpha, \mathcal{N}_0, \mathcal{N}_1, \ldots, \mathcal{N}_{k-1}\right). \end{aligned}$$

Therefore, we have that

$$\begin{aligned} &P(\tilde{v}_j^{(k+1)} > \delta\sqrt{l}, \mathcal{N}_0, \ldots, \mathcal{N}_{k-1}) \\ &\quad \leq P(v_j(1) > \delta\sqrt{l}\,) \times P(\mathcal{N}_0, \ldots, \mathcal{N}_{k-1}) \\ &\quad \leq \frac{E|v_j(1)|^{2+2\varepsilon_1}}{\delta^{2+2\varepsilon_1}} \frac{1}{l^{1+\varepsilon_1}} P(\mathcal{N}_0, \ldots, \mathcal{N}_{k-1}). \end{aligned}$$

Combining these, we conclude

$$\begin{aligned} &P(D_k^c, \mathcal{N}_{k-1}, \ldots, \mathcal{N}_0) \\ &\quad \leq \frac{P(\mathcal{N}_0, \ldots, \mathcal{N}_{k-1})}{\delta^{2+2\varepsilon_1}} \left[\sum_{i \in \mathcal{A}} E|u_i(1)|^{2+2\varepsilon_1} + \sum_{j=1}^{n} E|v_j(1)|^{2+2\varepsilon_1}\right] \frac{1}{l^{1+\varepsilon_1}}. \end{aligned} \tag{121}$$

We now consider $P(C_k, \mathcal{N}_{k-1}, \ldots, \mathcal{N}_0)$:

$$\begin{aligned} &P(C_k, \mathcal{N}_{k-1}, \ldots, \mathcal{N}_0) \\ &\quad = P\Bigg(\mathcal{N}_{k-1}, \ldots, \mathcal{N}_0, \\ &\qquad\qquad S_j\left(\sum_{i=0}^{k} x_j^{(i)} l\right) - S_j\left(\sum_{i=0}^{k-1} x_j^{(i)} l\right) \geq (2\mu_j + 1)l \text{ for some } j\Bigg) \end{aligned}$$

$$\leq \sum_{j=1}^{n} P\left(\tilde{v}_j^{(k)} + \sum_{i=S_j(\sum_{i=0}^{k-1} x_j^{(i)} l)+2}^{S_j(\sum_{i=0}^{k-1} x_j^{(i)} l)+\lceil (2\mu_j+1)l\rceil} v_j(i) \leq x_j^{(k)} l, \mathcal{N}_{k-1}, \ldots, \mathcal{N}_0\right)$$

$$\leq \sum_{j=1}^{n} \sum_{\alpha=0}^{\infty} P\left(\sum_{i=\alpha+2}^{\alpha+\lceil (2\mu_j+1)l\rceil} v_j(i) \leq l, S_j\left(\sum_{i=0}^{k-1} x_j^{(i)} l\right) = \alpha, \mathcal{N}_{k-1}, \ldots, \mathcal{N}_0\right)$$

$$\leq \sum_{j=1}^{n} \sum_{\alpha=0}^{\infty} \left[P\left(\sum_{i=\alpha+2}^{\alpha+\lceil (2\mu_j+1)l\rceil} v_j(i) \leq l\right) \times P\left(S_j\left(\sum_{i=0}^{k-1} x_j^{(i)} l\right) = \alpha, \mathcal{N}_{k-1}, \ldots, \mathcal{N}_0\right)\right]$$

$$\leq P(\mathcal{N}_{k-1}, \ldots, \mathcal{N}_0) \sum_{j=1}^{n} P(V_j(\lfloor (2\mu_j+1)l \rfloor) \leq l)$$

$$\leq \sum_{j=1}^{n} P(|S(el) - \mu l| \geq l) \qquad \text{for } l > \max_{j=1,\ldots,n}\{1/\mu_j\}.$$

Therefore, we conclude by Lemma 9 that

$$P(C_k, \mathcal{N}_{k-1}, \ldots, \mathcal{N}_0) \leq \frac{nC_3(1)}{l^{1+\varepsilon_1}} P(\mathcal{N}_{k-1}, \ldots, \mathcal{N}_0) \qquad \text{for } l > 2 \vee \max_{j=1,\ldots,n}\{1/\mu_j\}. \tag{122}$$

Next, we consider $P(B_k^c, C_k^c, \mathcal{N}_{k-1}, \ldots, \mathcal{N}_0)$:

$$\begin{aligned} &P(B_k^c, C_k^c, \mathcal{N}_{k-1}, \ldots, \mathcal{N}_0) \\ &\quad \leq \sum_{i=2}^{m} P\left(\sup_{\tau_k \leq s \leq \tau_{k+1}} Z_i(s) > [2|\theta^*| + n(2|\mu| + 1) + 2|\lambda|(1 + y'\theta^*)]l,\right. \\ &\qquad\qquad \left. \mathcal{N}_0, \ldots, \mathcal{N}_{k-1}, S_j\left(\sum_{i=1}^{k} x_j^{(i)} l\right) - S_j\left(\sum_{i=1}^{k} x_j^{(i)} l\right) \leq (2\mu_j + 1)l \ \forall j\right). \end{aligned}$$

For $\tau_k \leq s \leq \tau_{k+1}$, $i = 2, \ldots, m$, we have on the set $\mathcal{N}_0 \cap \cdots \cap \mathcal{N}_{k-1} \cap C_k^c$ that

$$\begin{aligned} Z_i(s) &\leq Z_i(\tau_1) + E_i(\tau_{k+1}) - E_i(\tau_k) + \sum_{j=1}^{n} \Phi_{ji}\left(S_j\left(\sum_{i=0}^{k} x_j^{(i)} l\right) - S_j\left(\sum_{i=0}^{k-1} x_j^{(i)} l\right)\right) \\ &\leq 2\theta_i + n(2|\mu| + 1)l + E_i(\tau_{k+1}) - E_i(\tau_k). \end{aligned}$$

Therefore,

$$
\begin{aligned}
&P(B_k^c, C_k^c, \mathcal{N}_{k-1}, \ldots, \mathcal{N}_0) \\
&\quad \le \sum_{i=2}^m P(E_i(\tau_{k+1}) - E_i(\tau_k) \ge 2|\lambda|(1 + y'\theta^*)l, \mathcal{N}_0, \ldots, \mathcal{N}_{k-1}) \\
&\quad \le \sum_{i=2}^m P(E_i(\tau_k + l(1 + y'\theta^*)) - E_i(\tau_k) \ge 2|\lambda|(1 + y'\theta^*)l, \mathcal{N}_0, \ldots, \mathcal{N}_{k-1}) \\
&\quad \le \sum_{i=2}^m P\Bigg(\tilde{u}_i^{(k)} + \sum_{\alpha=E_i(\tau_k)+2}^{E_i(\tau_k)+\lceil 2|\lambda|l(1+y'\theta^*)\rceil} u_i(\alpha) \le (1 + y'\theta^*)l, \mathcal{N}_0, \ldots, \mathcal{N}_{k-1}\Bigg) \\
&\quad \le \sum_{i=2}^m P\Bigg(\sum_{\alpha=E_i(\tau_k)+2}^{E_i(\tau_k)+\lceil 2|\lambda|l(1+y'\theta^*)\rceil} u_i(\alpha) \le (1 + y'\theta^*)l, \mathcal{N}_0, \ldots, \mathcal{N}_{k-1}\Bigg) \\
&\quad \le P(\mathcal{N}_0, \ldots, \mathcal{N}_{k-1}) \sum_{i=2}^m P(U_i(\lfloor 2|\lambda|(1 + y'\theta^*)l \rfloor) \le (1 + y'\theta^*)l) \\
&\quad \le P(\mathcal{N}_0, \ldots, \mathcal{N}_{k-1}) \sum_{i=2}^m P\bigg(|E(1 + y'\theta^*)l - \lambda(1 + y'\theta^*)l| \ge \frac{|\lambda|(1 + y'\theta^*)}{2} l\bigg),
\end{aligned}
$$

for $l > 1/|\lambda|$. We then have by Lemma 9 that

$$
\begin{aligned}
&P(B_k^c, C_k^c, \mathcal{N}_{k-1}, \ldots, \mathcal{N}_0) \\
&\quad \le \frac{mC_2((|\lambda|(1 + y'\theta^*))/2)}{(1 + y'\theta^*)^{1+\varepsilon_1}} P(\mathcal{N}_0, \ldots, \mathcal{N}_{k-1}) \frac{1}{l^{1+\varepsilon_1}}, \\
&\qquad\qquad l > \frac{1}{|\lambda|}\bigg[\frac{4}{1 + y'\theta^*} \vee 1\bigg].
\end{aligned}
\tag{123}
$$

Next, we consider $P(E_k^c, C_k^c, \mathcal{N}_{k-1}, \ldots, \mathcal{N}_0)$. On the set $C_k^c \cap \mathcal{N}_{k-1} \cap \cdots \cap \mathcal{N}_0$, one has that

$$
\int_{\tau_k}^{\tau_{k+1}} \mathbb{1}_{\{W(s) > (y'\theta^* + mn|y|(2|\mu|+1) + 2|y|my'\theta^*|\lambda|)l\}}\, dI_W(s) = I_1 + I_2,
$$

where

$$
I_1 = \int_{\tau_k}^{\tau_k + [y'(\theta - Z(\tau_k))]^+} \mathbb{1}_{\{W(s) > (y'\theta^* + mn|y|(2|\mu|+1) + 2|y|my'\theta^*|\lambda|)l\}}\, dI_W(s)
$$

and

$$
I_2 = \int_{\tau_k + [y'(\theta - Z(\tau_k))]^+}^{\tau_{k+1}} \mathbb{1}_{\{W(s) > (y'\theta^* + mn|y|(2|\mu|+1) + 2|y|my'\theta^*|\lambda|)l\}}\, dI_W(s).
$$

It is essential to observe that $I_2 = 0$ on the set $C_k^c \cap \mathcal{N}_{k-1} \cap \cdots \cap \mathcal{N}_0$, which is a consequence of the fact that $I_W(\tau_k + [y'(\theta - Z(\tau_k))]^+) = I_W(\tau_{k+1})$ on that set. The latter assertion follows because the servers work continuously during $[\tau_k + [y'(\theta - Z(\tau_k))]^+, \tau_{k+1}]$ by our policy description. (As we restrict attention on the set $C_k^c \cap \mathcal{N}_{k-1} \cap \cdots \cap \mathcal{N}_0$, the servers will have enough input to work on during $[\tau_k + [y'(\theta - Z(\tau_k))]^+, \tau_{k+1}]$.) Notice also that

$$\mathbb{1}_{\{I_1>0\}} \le \mathbb{1}_{\{\sup_{\tau_k \le s \le \tau_k + [y'(\theta - Z(\tau_k))]^+} W(s) > (y'\theta^* + mn|y|(2|\mu|+1) + 2|y|my'\theta^*|\lambda|)l\}}. \tag{124}$$

Without loss of generality, we can assume $y'\theta > y'Z(\tau_k) = W(\tau_k)$. (Otherwise, $I_1 = 0$.) For $s \in [\tau_k, \tau_k + [y'(\theta - Z(\tau_k))]^+]$, it is easy to see that

$$Z_i(s) - Z_i(\tau_k) \le E_i(s) - E_i(\tau_k) + \sum_{j=1}^{n} \Phi_{ji}\left(S_j\left(\sum_{i=0}^{k} x_j^{(i)} l\right) - S_j\left(\sum_{i=0}^{k-1} x_j^{(i)} l\right)\right), \qquad i = 1, \ldots, m.$$

Then it is straightforward to arrive at the following:

$$\begin{aligned} W(s) &\le W(\tau_k) + |y| \sum_{i=1}^{m} [E_i(s) - E_i(\tau_k)] + mn|y|(2|\mu| + 1)l \\ &\le y'\theta^* l + mn(2|\mu| + 1)l + |y| \sum_{i=1}^{m} [E_i(s) - E_i(\tau_k)]. \end{aligned}$$

Therefore, we immediately have the following:

$$\begin{aligned} &P\left(\sup_{\tau_k \le s \le \tau_k + [y'(\theta - Z(\tau_k))]^+} W(s) > (y'\theta^* + mn|y|(2|\mu| + 1) + 2|y|my'\theta^*|\lambda|)l\right) \\ &\quad \le P\left(\sum_{i=1}^{m} [E_i(\tau_k + y'\theta^* l) - E_i(\tau_k)] \ge 2mly'\theta^*|\lambda|\right) \\ &\quad \le \sum_{i \in \mathcal{A}} P([E_i(\tau_k + y'\theta^* l) - E_i(\tau_k)] \ge 2mly'\theta^*|\lambda|) \\ &\quad = \sum_{i \in \mathcal{A}} P\left(\tilde{u}_i^{(k)} + \sum_{\alpha = E_i(\tau_k)+2}^{E_i(\tau_k + \lceil 2ly'\theta^*|\lambda| \rceil)} u_i(\alpha) \le y'\theta^* l\right) \\ &\quad \le \sum_{i \in \mathcal{A}} P\left(\sum_{\alpha = E_i(\tau_k)+2}^{E_i(\tau_k + \lceil 2ly'\theta^*|\lambda| \rceil)} u_i(\alpha) \le y'\theta^* l\right) \\ &\quad = \sum_{i \in \mathcal{A}} P(U_i(\lfloor 2ly'\theta^*|\lambda| \rfloor) \le y'\theta^* l) \\ &\quad \le \sum_{i \in \mathcal{A}} P(E_i(y'\theta^* l) \ge 2ly'\theta^*|\lambda| - 1) \\ &\quad \le \sum_{i \in \mathcal{A}} P\left(|E_i(y'\theta^* l) - y'\theta^* l\lambda_i| \ge \frac{ly'\theta^*\lambda_i}{2}\right) \end{aligned}$$

$$\leq \left[\sum_{i\in\mathcal{A}} C_2\left(\frac{\lambda_i}{2}\right)\right]\frac{1}{(y'^{\theta^*})^{1+\varepsilon_1}}\frac{1}{l^{1+\varepsilon_1}} \qquad \text{for } l > \max_{i\in\mathcal{A}}\left\{\frac{4}{y'\theta^*\lambda_i}\right\},$$

where the last inequality follows from Lemma 9. Therefore, it is straightforward to conclude by (124) that

$$P(E_k^c, C_k^c, \mathcal{N}_{k-1}, \dots, \mathcal{N}_0) \leq \left[\sum_{i\in\mathcal{A}} C_2\left(\frac{\lambda_i}{2}\right)\right]\frac{1}{(y'^{\theta^*})^{1+\varepsilon_1}}\frac{1}{l^{1+\varepsilon_1}} \quad \text{for } l > \max_{i\in\mathcal{A}}\left\{\frac{4}{y'\theta^*\lambda_i}\right\}. \tag{125}$$

Finally, we consider $P(A_c^c, C_k^c, \mathcal{N}_{k-1}, \dots, \mathcal{N}_0)$:

$$P(A_k^c, C_k^c, \mathcal{N}_{k-1}, \dots, \mathcal{N}_0) \leq p_I^k + p_{II}^k + p_{III}^k, \tag{126}$$

where

$$p_I^k = P(|E(\tau_{k+1}) - E(\tau_k) - \lambda(\tau_{k+1} - \tau_k)| \geq \varepsilon l, \mathcal{N}_{k-1}, \dots, \mathcal{N}_0),$$

$$p_{II}^k = P\left(\sum_{j=1}^n \left|\Phi_j\left(S_j\left(\sum_{i=0}^k x_j^{(i)} l\right) - S_j\left(\sum_{i=0}^{k-1} x_j^{(i)} l\right)\right)\right.\right.$$
$$\left.\left. - P_j'\left(S_j\left(\sum_{i=0}^k x_j^{(i)} l\right) - S_j\left(\sum_{i=0}^{k-1} x_j^{(i)} l\right)\right)\right| \geq \varepsilon l, \mathcal{N}_{k-1}, \dots, \mathcal{N}_0, C_k^c\right),$$

$$p_{III}^k = \left(|RM|\left|S_j\left(\sum_{i=0}^k x_j^{(i)} l\right) - S_j\left(\sum_{i=0}^{k-1} x_j^{(i)} l\right) - M^{-1}x^{(k)} l\right| \geq \varepsilon l, \mathcal{N}_{k-1}, \dots, \mathcal{N}_0\right).$$

We first consider $p_{II}^k$, and immediately write

$$p_{II}^k \leq \sum_{j=1}^n P\left(\sup_{i=1,\dots,\lceil(2\mu_j+1)l\rceil} |\Phi_j(i) - iP_j'| \geq \frac{\varepsilon l}{n}, \mathcal{N}_0, \mathcal{N}_1, \dots, \mathcal{N}_{k-1}\right).$$

Therefore, we have by Lemma 10 that

$$p_{II}^k \leq \left[\sum_{j=1}^n \frac{C_{4j}(\varepsilon/(n(2\mu_j+2)))}{(2\mu_j+1)^{1+\varepsilon_1}}\right]\frac{1}{l^{1+\varepsilon_1}}. \tag{127}$$

We now consider $p_I^k$:

$$p_I^k \leq \sum_{i\in\mathcal{A}} [P_{I,a}^k(i) + P_{I,b}^k(i)], \tag{128}$$

where

$$P_{I,a}^k(i) = P(E_i(\tau_{k+1}) - E_i(\tau_k) \leq \lambda_i(\tau_{k+1} - \tau_k) - \varepsilon l, \mathcal{N}_0, \dots, \mathcal{N}_{k-1}),$$

$$P_{I,b}^k(i) = P(E_i(\tau_{k+1}) - E_i(\tau_k) \geq \lambda_i(\tau_{k+1} - \tau_k) + \varepsilon l, \mathcal{N}_0, \dots, \mathcal{N}_{k-1}).$$

We analyze each of these terms separately below:

$$P^k_{I,a}(i) = P\biggl(\tilde{u}^{(k)}_i + \sum_{j=E_i(\tau_k)+2}^{E_i(\tau_k)+\lfloor\lambda_i(\tau_{k+1}-\tau_k)-\varepsilon l\rfloor} u_i(j) \geq (\tau_{k+1}-\tau_k), \mathcal{N}_0,\ldots,\mathcal{N}_{k-1}\biggr)$$
$$\leq P\biggl(\delta\sqrt{l} + \sum_{j=E_i(\tau_k)+2}^{E_i(\tau_k)+\lfloor\lambda_i(\tau_{k+1}-\tau_k)-\varepsilon l\rfloor} u_i(j) \geq (\tau_{k+1}-\tau_k), \mathcal{N}_0,\ldots,\mathcal{N}_{k-1}\biggr).$$

By partitioning the event in question, first with respect to the values taken by $Z(\tau_k)$ and then with respect to the values taken by $E_i(\tau_k)$, we write

$$P^k_{I,a}(i) \leq \sum_{q_k \in B_{\delta l}(\theta)} \sum_{\alpha=0}^{\infty} P\biggl(\delta\sqrt{l} + \sum_{j=\alpha+2}^{\alpha+\lfloor\lambda_i(l+[y'(\theta-q_k)]^+)-\varepsilon l\rfloor} u_i(j) \geq (l+[y'(\theta-q_k)]^+), \mathcal{N}_0,\ldots,\mathcal{N}_{k-1}, Z(\tau_k)=q_k, E_i(\tau_k)=\alpha\biggr)$$
$$= \sum_{q_k \in B_{\delta l}(\theta)} \sum_{\alpha=0}^{\infty} P(U_i(\lfloor\lambda_i(l+[y'(\theta-q_k)]^+ - \varepsilon l\rfloor) \geq (l+[y'(\theta-q_k)]^+ - \delta\sqrt{l}\,)) \times P(\mathcal{N}_0,\ldots,\mathcal{N}_{k-1}, Z(\tau_k)=q_k, E_i(\tau_k)=\alpha)$$
$$\leq \sum_{q_k \in B_{\delta l}(\theta)} P(E_i(l+[y'(\theta-q)]^+ - \delta\sqrt{l}\,) \leq \lambda_i(l+[y'(\theta-q_k)]^+ - \varepsilon l - 1) \times P(\mathcal{N}_0,\ldots,\mathcal{N}_{k-1}, Z(\tau_k)=q_k)$$
$$\leq \sum_{q_k \in B_{\delta l}(\theta)} P\biggl(\sup_{0\leq s\leq l+y'\theta} |E_i(s)-\lambda_i s| \geq \frac{\varepsilon l}{2}\biggr) \times P(\mathcal{N}_0,\ldots,\mathcal{N}_{k-1}, Z(\tau_k)=q_k) \qquad \text{for } l > 144$$
$$\leq P\biggl(\sup_{0\leq s\leq l+y'\theta} |E_i(s)-\lambda_i s| \geq \frac{\varepsilon l}{2}\biggr) \sum_{q_k \in B_{\delta l}(\theta)} P(\mathcal{N}_0,\ldots,\mathcal{N}_{k-1}, Z(\tau_k)=q_k).$$

Hence, we have by Lemma 9 that

$$P^k_{I,a}(i) \leq \frac{C_2(\varepsilon/(2(1+y'\theta^*)))}{(1+y'\theta^*)^{1+\varepsilon_1}} P(\mathcal{N}_0,\ldots,\mathcal{N}_{k-1}) \frac{1}{l^{1+\varepsilon_1}} \tag{129}$$

$$\text{for } l > \frac{24}{\delta} \vee 144.$$

Next, we consider $P^k_{I,b}(i)$:

$$\begin{aligned} P^k_{I,b}(i) &= P(E_i(\tau_{k+1}) - E_i(\tau_k) \geq \lambda_i(\tau_{k+1} - \tau_k) + \varepsilon l, \mathcal{N}_0, \ldots, \mathcal{N}_{k-1}) \\ &\leq P\left(\tilde{u}_i^{(k)} + \sum_{j=E_i(\tau_k)+2}^{E_i(\tau_k)+\lceil\lambda_i(\tau_{k+1}-\tau_k)+\varepsilon l\rceil} u_i(j) \leq (\tau_{k+1} - \tau_k), \mathcal{N}_0, \ldots, \mathcal{N}_{k-1}\right) \\ &\leq P\left(\sum_{j=E_i(\tau_k)+2}^{E_i(\tau_k)+\lceil\lambda_i(\tau_{k+1}-\tau_k)+\varepsilon l\rceil} u_i(j) \leq (\tau_{k+1} - \tau_k), \mathcal{N}_0, \ldots, \mathcal{N}_{k-1}\right). \end{aligned}$$

It follows similarly that

$$P^k_{I,b}(i) \leq \frac{C_2(\varepsilon/(2(1+y'\theta^*)))}{(1+y'\theta^*)^{1+\varepsilon_1}} P(\mathcal{N}_0, \ldots, \mathcal{N}_{k-1}) \frac{1}{l^{1+\varepsilon_1}} \qquad \text{for } l > \frac{24}{\delta}. \tag{130}$$

By combining (128)–(130), we have

$$P^k_I \leq 2m \frac{C_2(\varepsilon/(2(1+y'\theta^*)))}{(1+y'\theta^*)^{1+\varepsilon_1}} P(\mathcal{N}_0, \ldots, \mathcal{N}_{k-1}) \frac{1}{l^{1+\varepsilon_1}} \qquad \text{for } l > \frac{24}{\delta} \vee 144. \tag{131}$$

To complete the proof we only need an estimate for $P^k_{III}(i)$. Observe that

$$P^k_{III} \leq \sum_{j=1}^{n} [P^k_{III,a}(j) + P^k_{III,b}(j)], \tag{132}$$

where

$$P^k_{III,a}(j) = P\left(S_j\left(\sum_{i=0}^{k} x_j^{(i)} l\right) - S_j\left(\sum_{i=0}^{k-1} x_j^{(i)} l\right) \geq \mu_j x_j^{(k)} l + \frac{\varepsilon l}{|RM|n}, \mathcal{N}_0, \ldots, \mathcal{N}_{k-1}\right),$$

$$P^k_{III,b}(j) = P\left(S_j\left(\sum_{i=0}^{k} x_j^{(i)} l\right) - S_j\left(\sum_{i=0}^{k-1} x_j^{(i)} l\right) \leq \mu_j x_j^{(k)} l - \frac{\varepsilon l}{|RM|n}, \mathcal{N}_0, \ldots, \mathcal{N}_{k-1}\right).$$

First consider $P^k_{III,a}(j)$:

$$
\begin{aligned}
P^k_{III,a}(j) &\le P\Bigg(\tilde v^{(k)}_j + \sum_{i=S_j(\sum_{i=0}^{k-1} x^{(i)}_j l)+2}^{S_j(\sum_{i=0}^{k-1} x^{(i)}_j l)+\lceil \mu_j x^{(k)}_j l+\varepsilon l/(|RM|n)\rceil} v_j(i) \le x^{(k)}_j l,\\
&\qquad\qquad \mathcal{N}_0,\dots,\mathcal{N}_{k-1}\Bigg)\\
&\le \sum_{\alpha=0}^{\infty} P\Bigg(\sum_{i=\alpha+2}^{\alpha+\lceil \mu_j x^{(k)}_j l+\varepsilon l/(|RM|n)\rceil} v_j(i) \le x^{(k)}_j l,\\
&\qquad\qquad \mathcal{N}_0,\dots,\mathcal{N}_{k-1}, S_j\Bigg(\sum_{i=0}^{k-1} x^{(i)}_j l\Bigg)=\alpha\Bigg)\\
&\le \sum_{\alpha=0}^{\infty} P\Bigg(S_j(x^{(k)}_j l) \ge \bigg\lfloor \mu_j x^{(k)}_j l + \frac{\varepsilon l}{|RM|n}\bigg\rfloor\Bigg)\\
&\qquad \times P\Bigg(\mathcal{N}_0,\dots,\mathcal{N}_{k-1}, S_j\Bigg(\sum_{i=0}^{k-1} x^{(i)}_j l\Bigg)=\alpha\Bigg)\\
&\le \sum_{\alpha=0}^{\infty} P\Bigg(\sup_{\|\tilde x\|\le 1,\tilde x\ge 0} |S(\tilde x l) - M^{-1}\tilde x l| \ge \frac{\varepsilon l}{2|RM|n}\Bigg)\\
&\qquad \times P\Bigg(\mathcal{N}_0,\dots,\mathcal{N}_{k-1}, S_j\Bigg(\sum_{i=0}^{k-1} x^{(i)}_j l\Bigg)=\alpha\Bigg).
\end{aligned}
$$

Therefore, we conclude by Lemma 9 that

$$
P^k_{III,a}(j) \le \frac{C_3(\varepsilon/(2|RM|n))}{l^{1+\varepsilon_1}} P(\mathcal{N}_0,\dots,\mathcal{N}_{k-1}) \qquad \text{for } l > \frac{24|RM|n}{\delta}. \tag{133}
$$

We now consider $P^k_{III,b}(j)$:

$$
\begin{aligned}
P^k_{III,b}(j) &= P\Bigg(\tilde v^{(k)}_j + \sum_{i=S_j(\sum_{i=0}^{k-1} x^{(i)}_j l)+2}^{S_j(\sum_{i=0}^{k-1} x^{(i)}_j l)+\lfloor \mu_j x^{(k)}_j l+\varepsilon l/(|RM|n)\rfloor} v_j(i) \ge x^{(k)}_j l,\\
&\qquad\qquad \mathcal{N}_0,\dots,\mathcal{N}_{k-1}\Bigg)\\
&\le \sum_{\alpha=0}^{\infty} P\Bigg(\sum_{i=\alpha+2}^{\alpha+\lfloor \mu_j x^{(k)}_j l+\varepsilon l/(|RM|n)\rfloor} v_j(i) \ge x^{(k)}_j l - \delta\sqrt{l},
\end{aligned}
$$

$$\mathcal{N}_0,\dots,\mathcal{N}_{k-1}, S_j\left(\sum_{i=0}^{k-1} x_j^{(i)} l\right) = \alpha\Bigg)$$

$$\leq \sum_{\alpha=0}^{\infty} P\left(S_j(x_j^{(k)} l - \delta\sqrt{l}\,) \leq \mu_j x_j^{(k)} l - \frac{\varepsilon l}{|RM|n}\right) \times P\left(\mathcal{N}_0,\dots,\mathcal{N}_{k-1}, S_j\left(\sum_{i=0}^{k-1} x_j^{(i)} l\right) = \alpha\right).$$

For $l > 144(|RM|n)^2$, we have that

$$P^k_{III,b}(j) \leq \sum_{\alpha=0}^{\infty} P\left(\sup_{\|\tilde{x}\|\leq 1, \tilde{x}\geq 0} |S(\tilde{x}l) - \mu\tilde{x}l| \geq \frac{\varepsilon l}{2|RM|n}\right) \times P\left(\mathcal{N}_0,\dots,\mathcal{N}_{k-1}, S_j\left(\sum_{i=0}^{k-1} x_j^{(i)} l\right) = \alpha\right).$$

Therefore, by Lemma 9 we conclude that

$$P^k_{III,b}(j) \leq \frac{C_3(\varepsilon/(2|RM|n))}{l^{1+\varepsilon_1}} P(\mathcal{N}_0,\dots,\mathcal{N}_{k-1}) \quad \text{for } l > \frac{24|RM|n}{\delta} \vee 144(|RM|n)^2. \tag{134}$$

By combining (132)–(134), we conclude

$$P^k_{III} \leq 2nC_2\left(\frac{\varepsilon}{2|RM|n}\right) l^{1+\varepsilon_1} P(\mathcal{N}_0,\dots,\mathcal{N}_{k-1}) \quad \text{for } l > \frac{24|RM|n}{\delta} \vee 144(|RM|n)^2. \tag{135}$$

Then by combining (127), (131) and (135), we conclude

$$P(A_k^c, C_k^c, \mathcal{N}_{k-1},\dots,\mathcal{N}_0) \leq \frac{C_A}{l^{1+\varepsilon_1}} P(\mathcal{N}_0,\dots,\mathcal{N}_{k-1}), \tag{136}$$

for $l > (24/\delta) \vee 144 \vee (24|RM|n/\delta) \vee 144(|RM|n^2)$, where

$$C_A = \frac{2m}{1+y'\theta^*} C_2\left(\frac{\varepsilon}{2(1+y'\theta^*)}\right) + 2nC_3\left(\frac{\varepsilon}{2|RM|n}\right) + \sum_{j=1}^{n} \frac{C_{4j}(\varepsilon/(n(2\mu_j+2)))}{(2\mu_j+1)^{1+\varepsilon_1}}.$$

Finally, by combining (121)–(123), (125) and (136), we conclude

$$P(\mathcal{N}_k^c, \mathcal{N}_0,\dots,\mathcal{N}_{k-1}) \leq \frac{C}{l^{1+\varepsilon_1}} P(\mathcal{N}_0,\dots,\mathcal{N}_{k-1}), \tag{137}$$

where

$$C = \widetilde{C} \vee l_{\min}, \tag{138}$$

and

$$\widetilde{C} \geq C_A + m\frac{C_2(|\lambda|(1+y'\theta^*)/2)}{(1+y'\theta^*)^{1+\varepsilon_1}} + nC_3(1) + \left[\sum_{i\in\mathcal{A}} C_2\left(\frac{\lambda_i}{2}\right)\right]\frac{1}{(y'\theta^*)^{1+\varepsilon_1}}$$
$$+ \frac{1}{\delta^{2+2\varepsilon_1}}\left[\sum_{i\in\mathcal{A}} E|u_i(1)|^{2+2\varepsilon_1} + \sum_{j=1}^{n} E|v_j(1)|^{2+2\varepsilon_1}\right],$$

and

$$l_{\min} = \frac{24}{\delta} \vee 144 \vee \frac{24|RM|n}{\delta} \vee 144(|RM|n)^2$$
$$\vee \max_{j=1,\dots,n}\left\{\frac{1}{\mu_j}\right\} \vee \frac{1}{|\lambda|}\left[\frac{4}{1+y'^{\theta^*}} \vee 1\right] \vee \max_{i\in\mathcal{A}}\left\{\frac{4}{y'\theta^*\lambda_i}\right\}. \qquad \square$$

A.4. *Proof of Proposition* 2. The following lemma will play an important role in proving Proposition 2.

LEMMA 11. *Given positive natural numbers* $J, M, L$, *let* $g$, $\{g_r\}_{r=1}^{\infty}$ *be functions from* $R^J$ *to* $R^M$, *where* $g$ *is continuous, such that*

$$g_r \to g \text{ uniformly over compact sets as } r \to \infty.$$

*Similarly, let* $f$, $\{f_r\}_{r=1}^{\infty}$ *be functions from* $R^M$ *to* $R^L$, *where* $f$ *is continuous, such that*

$$f_r \to f \text{ uniformly over compact sets as } r \to \infty.$$

*Then it follows that*

$$f_r \circ g_r \to f \circ g \text{ uniformly over compact sets as } r \to \infty.$$

PROOF. Let $K \subset R^J$ be an arbitrary compact set. We want to prove that

$$\lim_{r\to\infty} \sup_{x\in K} |f_r(g_r(x)) - f(g(x))| = 0.$$

Let $\varepsilon > 0$. First note that $g(K) = \{g(x) : x \in K \subset R^J\} \subset R^M$ is compact. This follows because $g$ is continuous, and $K$ is compact. Because $g(K)$ is compact, and $g_r \to g$ uniformly on $K$ as $r \to \infty$, there exists $R_1, N > 0$, such that

$$g_r(x) \in \widetilde{K} \qquad \forall x \in K, r \geq R_1 \qquad \text{where } \widetilde{K} = \{y \in R^M : |y| \leq N\}.$$

Then observe that for $r \geq R_1$, one has that

$$\sup_{x\in K} |f_r(g_r(x)) - f(g_r(x))| \leq \sup_{y\in \widetilde{K}} |f_r(y) - f(y)|.$$

Because $\widetilde{K} \subset R^M$ is compact, $f_r \to f$ uniformly on $\widetilde{K}$ as $r \to \infty$; and, therefore, there exists $R_2 \geq R_1$ such that

$$\sup_{x\in K} |f_r(g_r(x)) - f(g_r(x))| \leq \sup_{y\in \widetilde{K}} |f_r(y) - f(y)| < \frac{\varepsilon}{2} \qquad \text{for } r \geq R_2. \tag{139}$$

On the other hand, because $\widetilde{K}$ is compact, $f$ is uniformly continuous on $\widetilde{K}$. Therefore, there exists a $\delta > 0$ such that

$$|f(y_1) - f(y_2)| < \frac{\varepsilon}{2} \qquad \text{whenever } |y_1 - y_2| < \delta.$$

Moreover, because $g_r \to g$ uniformly on $K$ as $r \to \infty$, there exists $R_3 \geq R_2$ such that

$$\sup_{x\in K} |g_r(x) - g(x)| < \delta \qquad \text{for } r \geq R_3.$$

Therefore,

$$\sup_{x\in K} |f(g_r(x)) - f(g(x))| < \frac{\varepsilon}{2} \qquad \text{for } r \geq R_3. \tag{140}$$

Finally, observe that combining (139) and (140) gives the result. □

PROOF OF PROPOSITION 2. Fix $t > 0$ and $x > 0$. Let $\{T^r\}$ be an arbitrary sequence of admissible policies. Let $\{T^{r_j}\}$ be a subsequence such that

$$\lim_{j\to\infty} P(h \cdot \widehat{Z}_T^{r_j}(t) > x) = \liminf_{r\to\infty} P(h \cdot \widehat{Z}_T^{r}(t) > x).$$

Since $T^r$ corresponds to cumulative time allocations, $T^r$ is uniformly Lipschitz continuous with Lipschitz constant less than or equal to 1 for each $r$, and this property is preserved by the fluid scaled processes $\overline{T}^r(\cdot)$. It follows from this, and the functional central limit theorem for renewal processes (see [23]) and the fact that $\overline{T}_{DR} \Rightarrow x^*(\cdot)$ as $r \to \infty$, where $x^*(s) = x^* s$, $s \geq 0$ (and the fact that this limit is deterministic), that

$$\{(\widehat{E}^{r_j}(\cdot), \widehat{S}^{r_j}(\cdot), \widehat{\Phi}_1^{r_j}(\cdot), \ldots, \widehat{\Phi}_n^{r_j}(\cdot), \overline{T}^{r_j}(\cdot), \overline{T}_{DR}^{r_j}(\cdot))\}$$

is tight and any weak limit of this sequence has continuous paths almost surely (cf. Theorem 15.1 of [4]). In particular, the limit is of the following form:

$$(E^*(\cdot), S^*(\cdot), \Phi_1^*(\cdot), \ldots, \Phi_n^*(\cdot), \overline{T}(\cdot), x^*(\cdot)), \tag{141}$$

where $E^*(\cdot), S^*(\cdot), \Phi_1^*(\cdot), \ldots, \Phi_n^*(\cdot)$ are driftless Brownian motions of appropriate dimension, and $\overline{T}(\cdot)$ is a nondecreasing process with

$$A(\overline{T}(s_2) - \overline{T}(s_1)) \le (s_2 - s_1)e \qquad \text{for } 0 \le s_1 \le s_2 \text{ almost surely.} \tag{142}$$

Let $\{T^{r'}\}$ be a further subsequence of $\{T^{r_j}\}$ which converges weakly to a limit as in (141). By appealing to the Skorohod representation theorem, we may choose an equivalent distributional representation (which we will denote by putting a "˜" above the symbols) such that the sequence random processes

$$\{(\tilde{\hat{E}}^{r'}(\cdot), \tilde{\hat{S}}^{r'}(\cdot), \tilde{\hat{\Phi}}_1^{r'}(\cdot), \ldots, \tilde{\hat{\Phi}}_n^{r'}(\cdot), \tilde{\overline{T}}^{r'}(\cdot), \tilde{\overline{T}}_{DR}^{r'}(\cdot))\},$$

as well as the limit

$$(\tilde{E}^*(\cdot), \tilde{S}^*(\cdot), \tilde{\Phi}_1^*(\cdot), \ldots, \tilde{\Phi}_n^*(\cdot), \tilde{\overline{T}}(\cdot), x^*(\cdot)),$$

are defined on a new probability space, say $(\tilde{\Omega}, \tilde{\mathcal{F}}, \tilde{P})$, so that $\tilde{P}$ almost surely

$$\begin{aligned}(\tilde{\hat{E}}^{r'}(\cdot), \tilde{\hat{S}}^{r'}(\cdot), \tilde{\hat{\Phi}}_1^{r'}(\cdot), \ldots, \tilde{\hat{\Phi}}_n^{r'}(\cdot), \tilde{\overline{T}}^{r'}(\cdot), \tilde{\overline{T}}_{DR}^{r'}(\cdot))& \\ \to (\tilde{E}^*(\cdot), \tilde{S}^*(\cdot), \tilde{\Phi}_1^*(\cdot), \ldots, \tilde{\Phi}_n^*(\cdot), \tilde{\overline{T}}(\cdot), x^*(\cdot))&\end{aligned} \tag{143}$$

uniformly over compact time intervals as $r' \to \infty$.

Remark. Two processes are said to be equal in distribution if they have the same finite-dimensional distributions. This is also referred to as one is a version of the other; see page 50 of [11].

We define the following processes on this new probability space:

$$\tilde{\overline{E}}^{r'}(s) = \frac{1}{r'}\tilde{\hat{E}}^{r'}(s) + \lambda s, \qquad s \ge 0, \tag{144}$$

$$\tilde{\overline{S}}^{r'}(s) = \frac{1}{r'}\tilde{\hat{S}}^{r'}(s) + \mu s, \qquad s \ge 0, \tag{145}$$

$$\tilde{\overline{\Phi}}_j^{r'}(s) = \frac{1}{r'}\tilde{\hat{\Phi}}_j^{r'}(s) + P_j' s, \qquad s \ge 0, j = 1, \ldots, n, \tag{146}$$

$$\tilde{\overline{Z}}^{r'}(0) = \frac{\theta(r')}{r'}. \tag{147}$$

We note that these processes have the same joint distribution as the corresponding processes in the old probability space. We also define the following processes associated with the sequence of (scaled) policies $\{\tilde{\overline{T}}^{r'}(\cdot)\}$ on the

new probability space:

$$\tilde{\hat{T}}^{r'}(s) = r'\tilde{\bar{T}}^{r'}(s), \qquad s \geq 0, \tag{148}$$

$$\tilde{\hat{Y}}^{r'}_T(s) = r'x^*s - \tilde{\hat{T}}^{r'}(s), \qquad s \geq 0, \tag{149}$$

$$\tilde{\hat{I}}^{r'}_T(s) = r'es - A\tilde{\hat{T}}^{r'}(s), \qquad s \geq 0, \tag{150}$$

$$\begin{aligned} \tilde{\hat{X}}^{r'}_T(s) &= \tilde{\hat{Z}}^{r'}(0) + \tilde{\hat{E}}^{r'}(s) \\ &\quad + \sum_{j=1}^{n} \tilde{\hat{\Phi}}^{r'}_j(\tilde{\bar{S}}^{r'}(\tilde{\bar{T}}^{r'}(s))) - RM\tilde{\hat{S}}^{r'}(\tilde{\bar{T}}^{r'}(s)), \qquad s \geq 0, \end{aligned} \tag{151}$$

$$\tilde{\hat{Z}}^{r'}_T(s) = \tilde{\hat{X}}^{r'}_T(s) + R\tilde{\hat{Y}}^{r'}_T(s), \qquad s \geq 0, \tag{152}$$

$$\tilde{\hat{I}}^{r'}_{W,T}(s) = \pi'\tilde{\hat{I}}^{r'}_T(s), \qquad s \geq 0, \tag{153}$$

$$\tilde{\hat{I}}^{r'}_{N,T}(s) = -\eta'\tilde{\hat{Y}}^{r'}_{T,N}(s), \qquad s \geq 0, \tag{154}$$

where $\tilde{\hat{Y}}^{r'}_{T,N}$ is the last $n-b$ components of $\tilde{\hat{Y}}^{r'}_T$. Furthermore,

$$\tilde{\hat{X}}^{r'}_{W,T}(s) = y'\tilde{\hat{X}}^{r'}_T(s), \qquad s \geq 0, \tag{155}$$

$$\widetilde{\widehat{W}}^{r'}_T(s) = \tilde{\hat{X}}^{r'}_{W,T}(s) + \tilde{\hat{I}}^{r'}_{W,T}(s) + \tilde{\hat{I}}^{r'}_{W,N}(s), \qquad s \geq 0, \tag{156}$$

$$\tilde{\bar{Z}}^{r'}_T(s) = \frac{1}{r'}\tilde{\hat{Z}}^{r'}_T(s), \qquad s \geq 0, \tag{157}$$

$$\widetilde{\overline{W}}^{r'}_T(s) = \frac{1}{r'}\widetilde{\widehat{W}}^{r'}_T(s), \qquad s \geq 0. \tag{158}$$

It follows from these, Lemma 1 and definition of $\eta$ (by purely algebraic and straightforward manipulations) that $\widetilde{\widehat{W}}^{r'}_T(\cdot)$ has the following equivalent representation:

$$\widetilde{\widehat{W}}^{r'}_T(s) = y'\tilde{\hat{Z}}^{r'}_T(s), \qquad s \geq 0. \tag{159}$$

It is crucial to observe that these processes defined by (144)–(158) have the same joint distribution as the corresponding scaled processes in the old probability space for each $r$.

Because in the old space we have that almost surely $\overline{T}^{r'}(\cdot)$ is nondecreasing, and $\tilde{\bar{T}}^{r'}(\cdot)$ and $\overline{T}^{r'}(\cdot)$ are equal in distribution in $D^n[0,\infty)$, we have that

(in the new probability space)

$$\tag{160} \tilde{\bar{T}}^{r'}(\cdot) \text{ is nondecreasing almost surely.}$$

Similarly, we conclude that almost surely (in the new probability space)

$$\tag{161} A(\tilde{\bar{T}}^{r'}(s_2) - \tilde{\bar{T}}^{r'}(s_1)) \le (s_2 - s_1)e, \qquad 0 \le s_1 \le s_2,$$

$$\tag{162} \tilde{\hat{T}}^{r'}(\cdot) \text{ is nondecreasing,}$$

$$\tag{163} \tilde{\hat{I}}_T^{r'}(\cdot) \text{ is nondecreasing.}$$

We also define the sequence of processes $\{\tilde{\hat{X}}^{r'}_{DR}(\cdot)\}$ associated with our discrete review policy as follows:

$$\tag{164} \begin{aligned} \tilde{\hat{X}}^{r'}_{DR}(s) &= \tilde{\hat{Z}}^{r'}(0) + \tilde{\hat{E}}^{r'}(s) \\ &\quad + \sum_{j=1}^{n} \tilde{\hat{\Phi}}_j^{r'}(\tilde{\bar{S}}^{r'}(\tilde{\bar{T}}^{r'}_{DR}(s))) - RM\tilde{\hat{S}}^{r'}(\tilde{\bar{T}}^{r'}_{DR}(s)), \qquad s \ge 0. \end{aligned}$$

[It is essential to observe that the processes defined by (144)–(158) and (164), as well as the processes $(\tilde{\hat{E}}^{r'}(\cdot), \tilde{\hat{S}}^{r'}(\cdot), \tilde{\hat{\Phi}}_1^{r'}(\cdot), \ldots, \tilde{\hat{\Phi}}_n^{r'}(\cdot), \tilde{\bar{T}}^{r'}(\cdot), \tilde{\bar{T}}^{r'}_{DR}(\cdot))$, have the same joint distribution as the corresponding processes defined on the old probability space.]

First, because $\hat{Z}_T^{r'}(\cdot)$ is a nonnegative process almost surely on the old space, and $\hat{Z}_T^{r'}(\cdot)$ and $\tilde{\hat{Z}}_T^{r'}(\cdot)$ are equal in distribution, it follows that $\tilde{\hat{Z}}_T^{r'}(\cdot)$ is nonnegative almost surely. Similarly, we argue that $\widetilde{\widehat{W}}_T^{r'}(\cdot)$ is nonnegative almost surely too. It also follows from (149), (153), (154), (162) and (163) that $\tilde{\hat{I}}^{r'}_{W,T}(\cdot)$ and $\tilde{\hat{I}}^{r'}_{W,N}(\cdot)$ are nondecreasing processes. Therefore, it follows from these, (156) and (159) that

$$\tag{165} \begin{aligned} h \cdot \tilde{\hat{Z}}_T^{r'}(t) &\ge \frac{h_1}{y_1} \widetilde{\widehat{W}}_T^{r'}(t) \\ &\ge \frac{h_1}{y_1}(\tilde{\hat{X}}^{r'}_{W,T}(t) + \tilde{\hat{I}}^{r'}_{W,T}(t) + \tilde{\hat{I}}^{r'}_{W,N}(t)) \\ &\ge \frac{h_1}{y_1}\varphi(\tilde{\hat{X}}^{r'}_{W,T})(t), \end{aligned}$$

where the last statement follows because $\widetilde{\widehat{W}}_T^{r'}(\cdot)$ is a nonnegative process, and $\tilde{\hat{I}}^{r'}_{W,T}(\cdot)$ and $\tilde{\hat{I}}^{r'}_{W,N}(\cdot)$ are nondecreasing processes; see Appendix B of [1].

From (144) we immediately conclude that $\widetilde{P}$ almost surely as $r' \to \infty$,

$$(\widetilde{\widetilde{S}}^{r'}(\cdot), \widetilde{\widetilde{E}}^{r'}(\cdot), \widetilde{\widetilde{\Phi}}_1^{r'}(\cdot), \ldots, \widetilde{\widetilde{\Phi}}_n^{r'}(\cdot)) \to (\mu(\cdot), \lambda(\cdot), P_1(\cdot), \ldots, P_n(\cdot)) \qquad \text{u.o.c.,} \tag{166}$$

where

$$\lambda(s) = \lambda s, \qquad s \geq 0,$$
$$\mu(s) = \mu s, \qquad s \geq 0,$$
$$P_j(s) = P_j' s, \qquad s \geq 0, j = 1, \ldots, n,$$

and u.o.c. means uniformly over compact time intervals. Furthermore,

$$\widetilde{\widetilde{Z}}_T^{r'}(\cdot) \to \widetilde{\widetilde{Z}}_T(\cdot) = \lambda(\cdot) - R\widetilde{\widetilde{T}}(\cdot) \qquad \text{u.o.c.,} \tag{167}$$

$$\widetilde{\widetilde{W}}_T^{r'}(\cdot) \to \widetilde{\widetilde{W}}_T(\cdot) = y'\widetilde{\widetilde{Z}}_T(\cdot) \qquad \text{u.o.c.,} \tag{168}$$

and it follows by Lemma 11 that

$$\widetilde{\widetilde{X}}_T^{r'} \to \widetilde{X}(\cdot) = \widetilde{E}^*(\cdot) - \sum_{j=1}^n \widetilde{\Phi}_j^*(\mu(\widetilde{\widetilde{T}}(\cdot))) - RM\widetilde{S}^*(\widetilde{\widetilde{T}}(\cdot)) \qquad \text{u.o.c.,} \tag{169}$$

$$\widetilde{\widetilde{X}}_{DR}^{r'} \to \widetilde{X}^*(\cdot) = \widetilde{E}^*(\cdot) - \sum_{j=1}^n \widetilde{\Phi}_j^*(\mu(x^*(\cdot))) - RM\widetilde{S}^*(x^*(\cdot)) \qquad \text{u.o.c.,} \tag{170}$$

and

$$\widetilde{\widetilde{X}}_{W,T}^{r'}(\cdot) \to \widetilde{X}_W(\cdot) = y'\widetilde{X}(\cdot) \qquad \text{u.o.c.,} \tag{171}$$

$$\widetilde{\widetilde{X}}_{W,DR}^{r'}(\cdot) \to \widetilde{X}_{W,DR}(\cdot) = y'\widetilde{X}^*(\cdot) \qquad \text{u.o.c.} \tag{172}$$

The following conclusions are immediate from (161) and (162):

$$A(\widetilde{\widetilde{T}}(s_2) - \widetilde{\widetilde{T}}(s_1)) \leq (s_2 - s_1)e, \qquad 0 \leq s_1 \leq s_2, \tag{173}$$

$$\widetilde{\widetilde{T}}(\cdot) \text{ is nondecreasing.} \tag{174}$$

Let $\widetilde{\widetilde{T}}_N$ denote the last $n - b$ components of $\widetilde{\widetilde{T}}(\cdot)$. We expand (168) by using (167), and use (13), (16) and the definition of $\eta$ in (72) to obtain $\widetilde{\widetilde{W}}_T(s) = s - \pi' A\widetilde{\widetilde{T}}(s) + \eta\widetilde{\widetilde{T}}_N(s), s \geq 0$. Lemma 1 ensures that $\eta \geq 0$ and (174) ensures that $\widetilde{\widetilde{T}}_N(\cdot)$ is nondecreasing. Therefore, for any $t_2 \geq t_1 \geq 0$, one has that

$$\begin{aligned}\widetilde{\widetilde{W}}_T(t_2) - \widetilde{\widetilde{W}}_T(t_1) &= (t_2 - t_1) - \pi' A(\widetilde{\widetilde{T}}(t_2) - \widetilde{\widetilde{T}}(t_1)) + \eta(\widetilde{\widetilde{T}}_N(t_2) - \widetilde{\widetilde{T}}_N(t_1)) \\ &\geq (t_2 - t_1) - \pi' A(\widetilde{\widetilde{T}}(t_2) - \widetilde{\widetilde{T}}(t_1)).\end{aligned}$$

From (15), $\pi \geq 0$ and $\pi' e = 1$. Combining this with (173) yields

$$\widetilde{\overline{W}}_T(t_2) - \widetilde{\overline{W}}_T(t_1) \geq 0 \qquad \text{for all } t_2 \geq t_1 \geq 0. \tag{175}$$

That is, $\widetilde{\overline{W}}_T(\cdot)$ is nondecreasing.

Next, we fix $t > 0$ and define $U_t = \{\omega \in \widetilde{\Omega} : \widetilde{\overline{Z}}_T(t, \omega) = 0\}$. For $\widetilde{P}$ a.e. $\omega \in U_t$, we have that $\widetilde{\overline{W}}_T(t) = 0$. [This follows because $\widetilde{\overline{W}}_T(t) = y' \widetilde{\overline{Z}}_T(t)$, cf. (168).] It follows from this and (175) that $\widetilde{\overline{W}}_T(s) = 0, s \in [0, t]$. Combining this with the fact that $y' > 0$, we conclude that

$$\widetilde{\overline{Z}}_T(s) = 0, \qquad s \in [0, t].$$

Then we conclude by (167) that $\lambda s - R\widetilde{\overline{T}}(s) = 0, s \in [0, t]$. Combining this with (173), we see that $(\widetilde{\overline{T}}(s)/s, 1)$ is an optimal solution to the static planning problem by the heavy traffic assumption. By the uniqueness of this solution (also by the heavy traffic assumption), we see that $\widetilde{\overline{T}}(s) = x^* s, s \in [0, t]$. Combining this with (169) and (170), we have that

$$\widetilde{X}(s) = \widetilde{X}^*(s), \qquad s \in [0, t] \text{ for } \widetilde{P} \text{ a.e. } \omega \in U_t.$$

Also, by (171) and (172) we have that

$$\widetilde{X}_W(s) = \widetilde{X}^*_W(s), \qquad s \in [0, t] \text{ for } \widetilde{P} \text{ a.e. } \omega \in U_t.$$

Combining this with (165), (171) and continuity of the one-dimensional regulator map $\varphi(\cdot)$, we conclude that

$$\liminf_{r' \to \infty} h \cdot \widetilde{\widehat{Z}}^{r'}_T(t) \geq \frac{h_1}{y_1} \varphi(\widetilde{X}^*_W)(t) \qquad \text{for } \widetilde{P} \text{ a.e. } \omega \in U_t. \tag{176}$$

On the other hand, for $\widetilde{P}$ a.e. $\omega \in U^C_t$, we have that $\widetilde{\overline{Z}}_T(t) \neq 0$. Therefore, for $\widetilde{P}$ a.e. $\omega \in U^C_t$, we have $h \cdot \widetilde{\overline{Z}}_T(t) > 0$. Since $h \cdot \widetilde{\widehat{Z}}^{r'}_T(t) = r' h \cdot \widetilde{\overline{Z}}^{r'}_T(t)$ and $\lim_{r' \to \infty} r' h \cdot \widetilde{\overline{Z}}^{r'}_T(t) = \infty$ [the last statement follows because $\lim_{r' \to \infty} h \cdot \widetilde{\overline{Z}}^{r'}_T(t) = h \cdot \widetilde{\overline{Z}}_T(t) > 0$], we have that

$$\liminf_{r' \to \infty} h \cdot \widetilde{\widehat{Z}}^{r'}_T(t) = \infty \qquad \text{for } \widetilde{P} \text{ a.e. } \omega \in U^C_t.$$

Combining this with (176), we have that

$$\liminf_{r' \to \infty} h \cdot \widetilde{\widehat{Z}}^{r'}_T(t) \geq \frac{h_1}{y_1} \varphi(\widetilde{X}^*_W)(t) \qquad \text{for } \widetilde{P} \text{ a.e. } \omega \in \widetilde{\Omega}. \tag{177}$$

To conclude the proof we observe that

$$\begin{aligned}
\lim_{r'\to\infty} P(h\cdot \widehat{Z}_T^{r'}(t)>x) &= \lim_{r'\to\infty} \widetilde{P}(h\cdot \widetilde{\widehat{Z}}_T^{r'}(t)>x)\\
&= \lim_{r'\to\infty} \widetilde{E}[\mathbb{1}_{\{h\cdot \widetilde{\widehat{Z}}_T^{r'}(t)>x\}}]\\
&\geq \widetilde{E}\Big[\liminf_{r'\to\infty} \mathbb{1}_{\{h\cdot \widetilde{\widehat{Z}}_T^{r'}(t)>x\}}\Big]\\
&\geq \widetilde{E}[\mathbb{1}_{\{h_1/y_1\varphi(\widetilde{X}_W^*)(t)>x\}}]\\
&= \widetilde{P}\bigg(\frac{h_1}{y_1}\varphi(\widetilde{X}_W^*)(t)>x\bigg),
\end{aligned}$$

where the third step follows from Fatou's lemma. Moreover,

$$\widetilde{P}(\varphi(\widetilde{X}_W^*)(t)\geq x) = 2N\bigg(\frac{-xy_1}{h_1\sigma\sqrt{t}}\bigg),$$

where $N(\cdot)$ is the standard normal cumulative distribution function; see [13]. □

**Acknowledgments.** We are grateful to J. M. Harrison for sharing his ideas with us generously throughout the project and to R. J. Williams for a number of discussions on her related work. We also thank C. Maglaras for useful conversations during the initial phase of the project. Finally, we thank the Associate Editor and two anonymous referees for their careful reading and extensive comments which improved the paper considerably.

Kellogg School of Management
Northwestern University
Evanston, Illinois 60208-2009
USA
e-mail: b-ata@kellogg.northwestern.edu

Graduate School of Business
Stanford University
Stanford, California 94305-5015
USA
e-mail: skumar@stanford.edu